\documentclass[hidelinks,onefignum,onetabnum]{siamart251216}

\usepackage{amsfonts}
\usepackage{aliascnt} 

\headers{One-shot identification of analytic systems}{N. Boull\'e, D. Halikias, S. E. Otto, and A. Townsend}

\title{One-shot identification of analytic systems\thanks{Submitted \today.
    \funding{The work of A. T. was supported by the Defense Advanced Research Projects Agency (DARPA) through The Right Space (TRS) Disruption Opportunity (DARPA-PA-24-04-07) and National Science Foundation CAREER grant DMS-2045646.}}}

\author{Nicolas Boull\'e\thanks{Department of Mathematics, Imperial College London, London, SW7 2AZ, UK
    (\email{n.boulle@imperial.ac.uk}).}
  \and Diana Halikias\thanks{Courant Institute, New York University, New York, NY 10012, USA
    (\email{diana.halikias@nyu.edu}).}
  \and Samuel E. Otto\thanks{Sibley School of Mechanical and Aerospace Engineering, Cornell University, Ithaca, NY 14853-7501, USA
    (\email{s.otto@cornell.edu}).}
    \and Alex Townsend\thanks{Mathematics Department, Cornell University, Ithaca, NY 14853-4201, USA
    (\email{townsend@cornell.edu}).}}

\usepackage{amsopn}
\usepackage{mathtools}
\usepackage{booktabs}
\usepackage{bm}
\usepackage{enumitem}
\usepackage{algorithm}
\usepackage[noend]{algpseudocode}
\usepackage{cite}
\usepackage{comment}
\usepackage{tikz}
\usepackage{pgfplots}
\pgfplotsset{compat=1.18}
\usepackage{overpic}
\usepackage{standalone}
\usetikzlibrary{intersections,arrows.meta,positioning,calc,fit}
\usepackage{lineno}

\usepackage{jabbrv}
\DefineSpuriousJournalWord{on}
\DefineSpuriousJournalWord{for}
\DefineSpuriousJournalWord{Its}
\DefineSpuriousJournalWord{its}
\DefineJournalAbbreviation{Foundations}{Found}
\DefineJournalAbbreviation{Transactions}{Trans}

\newcommand{\R}{\mathbb{R}}
\newcommand{\C}{\mathbb{C}}

\newcommand{\St}{St}
\newcommand{\Fdeg}{\mcal{F}_{\rm{deg}}}
\newcommand{\Udeg}{\mcal{U}_{\rm{deg}}}
\renewcommand{\d}{\,\textup{d}}

\DeclareMathOperator*{\esssup}{ess\,sup}
\DeclareMathOperator{\AC}{AC}

\DeclareMathOperator{\rank}{rank}
\DeclareMathOperator{\Range}{Range}
\DeclareMathOperator{\Null}{Null}
\DeclareMathOperator{\vspan}{span}

\DeclareMathOperator{\jet}{J}
\DeclareMathOperator{\supp}{supp}
\DeclareMathOperator{\D}{D}

\newcommand{\vertiii}[1]{{\left\vert\kern-0.25ex\left\vert\kern-0.25ex\left\vert #1 \right\vert\kern-0.25ex\right\vert\kern-0.25ex\right\vert}}

\newcommand{\mcal}[1]{{\mathcal{#1}}}
\newcommand{\mat}[1]{{\bm{#1}}}
\newcommand{\vect}[1]{{\bm{#1}}}

\usepackage[dvipsnames]{xcolor}
\definecolor{c0}{RGB}{31, 119, 180}
\definecolor{c1}{RGB}{230, 180, 60}
\definecolor{c2}{RGB}{44, 160, 90}
\definecolor{c3}{RGB}{64, 170, 170}
\definecolor{c4}{RGB}{214, 85, 85}
\definecolor{c5}{RGB}{150, 90, 200}

\newaliascnt{assumption}{theorem}
\newtheorem{assumption}[assumption]{Assumption}
\aliascntresetthe{assumption}
\newaliascnt{example}{theorem}
\newtheorem{example}[example]{Example}
\aliascntresetthe{example}
\newaliascnt{remark}{theorem}
\newtheorem{remark}[remark]{Remark}
\aliascntresetthe{remark}

\crefformat{enumi}{#2#1#3}
\crefformat{enumii}{#2#1#3}
\crefformat{enumiii}{#2#1#3}
\crefformat{enumiv}{#2#1#3}

\crefname{assumption}{Assumption}{Assumptions}
\Crefname{assumption}{Assumption}{Assumptions}

\begin{document}

\maketitle

% REQUIRED
\begin{abstract}
  We study when the coefficients of a linearly parameterized analytic system can be identified from a single input-solution pair. Given a finite dictionary of analytic operators, each admissible input induces a finite-dimensional linear recovery map from coefficients to observed outputs. Coefficient recovery is possible exactly when this map is injective. We define the degenerate input set as the set of inputs for which injectivity fails and prove a zero-one law: either every admissible input is degenerate, or the degenerate set is negligible in the sense that it is compactly slice-null and has zero measure under every nondegenerate Gaussian measure. We then apply this framework to the recovery of matrices, analytic continuous-time dynamical systems, and differential operators. For several ODE and PDE classes, we further characterize the degenerate sets explicitly and identify deterministic choices of inputs that guarantee recovery.
\end{abstract}

% REQUIRED

% REQUIRED

\section{Introduction}
Differential equations provide a compact language for expressing physical laws, biological mechanisms, and engineered systems. Yet, for many complex systems, the governing physics remains only partially understood, motivating the identification of differential equations directly from data~\cite{karniadakis2021physics}. Data-driven model discovery has applications ranging from fluid dynamics~\cite{fukami2021sparse,brunton2020machine,vinuesa2022enhancing} and weather forecasting~\cite{pathak2022fourcastnet,lam2023learning,price2025probabilistic} to biological systems~\cite{mangan2017inferring}. A common approach is to posit a dictionary of candidate terms, such as derivatives, nonlinearities, or forcing terms, and then use observed input-solution pairs to determine which linear combination of these terms best explains the data.

In the general form of the problem, given a dictionary of analytic functions $\mathcal{D}_0,\mathcal{D}_1,\ldots,\mathcal{D}_N$, one seeks coefficients $(c_1^*,\ldots,c_N^*)\in\R^N$ such that the input-output data is governed by an equation:
\begin{equation} \label{eq:model}
  \sum_{i=1}^N c_i^* \mathcal{D}_i(u,f)=\mathcal{D}_0(u,f),
\end{equation}
where $c_i^*$ are unknowns and $u$ denotes the system state and $f$ denotes a forcing term or input. For example, a semilinear elliptic dictionary may contain $\mathcal{D}_1(u,f)=u_{xx}$, $\mathcal{D}_2(u,f)=u_{yy}$, and $\mathcal{D}_3(u,f)=u^3$; an ODE dictionary may contain integrated terms $\mathcal{D}_i[u,f](t)=\int_a^t D_i(s,u(s),f(s))\,\d s$; and a matrix-recovery dictionary may contain $\mathcal{D}_i(u,f)=A_i f$ for known matrices $A_i$. The goal is usually to recover the unknown coefficient vector $\vect{c}^*=(c_1^*,\ldots,c_N^*)$ from multiple observations of input-solution pairs $(f,u)$ satisfying \cref{eq:model}. This viewpoint underlies many modern approaches to data-driven model discovery and system identification, including symbolic regression and genetic programming~\cite{udrescu2020ai,udrescu2020ai2,schmidt2009distilling,searson2010gptips}, as well as algorithms that seek a sparse vector $\vect{c}$, such as Sparse Discovery of Nonlinear Dynamics (SINDy)~\cite{brunton2016discovering,champion2019data}. These methods often rely on sufficiently rich collections of trajectories, experiments, or time-series measurements to distinguish the active terms in the model. In many settings, however, one does not have access to large collections of trajectories or repeated experiments. Therefore, in this paper, we pose the fundamental one-shot identification question:
\medskip
\begin{quote}
  When does a \emph{single} $(f,u)$ pair determine the coefficients in~\cref{eq:model}?
\end{quote}
\medskip
Throughout, ``one-shot'' means one complete input-output experiment in the ambient model class. For a PDE, this may be one forcing term and its full solution field. For a dynamical system, it may be one forcing trajectory and the resulting state trajectory. For a matrix family, it may be one prescribed batch of matrix-vector queries. The question is whether a single chosen experiment induces an injective map from coefficients to observable outputs.

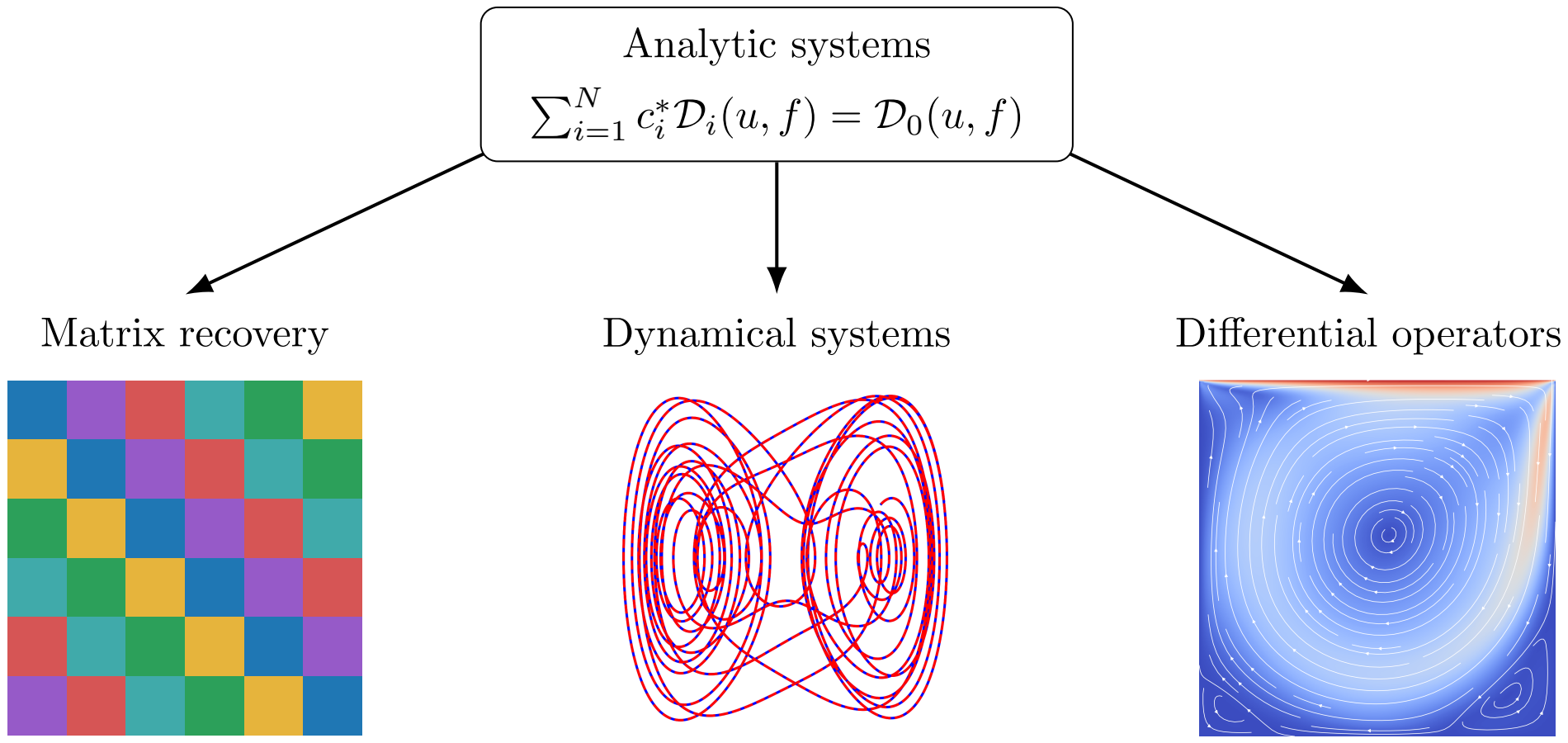
\begin{figure}[t]
  \centering
  \resizebox{\textwidth}{!}{
    \begin{tikzpicture}[node distance=1.5cm and 2.5cm,
    box/.style={rectangle, draw, rounded corners, align=center, minimum width=5cm, minimum height=1.3cm},
    arrow/.style={-{Latex}, thick},
    cloud/.style={draw, ellipse, minimum height=1cm, align=center}
    ]
    \pgfdeclarelayer{bg}
    \pgfsetlayers{bg,main}

    \node at (0,0) (matrix)
    {\includestandalone[height=3cm]{Figure/circulant}};
    \node at ($(matrix.north)+(0,0.25)$) {\hyperref[sec:matrix_recovery]{Matrix recovery}};

    \node at (5,0) (dynsys)
    {\begin{overpic}[height=3cm]{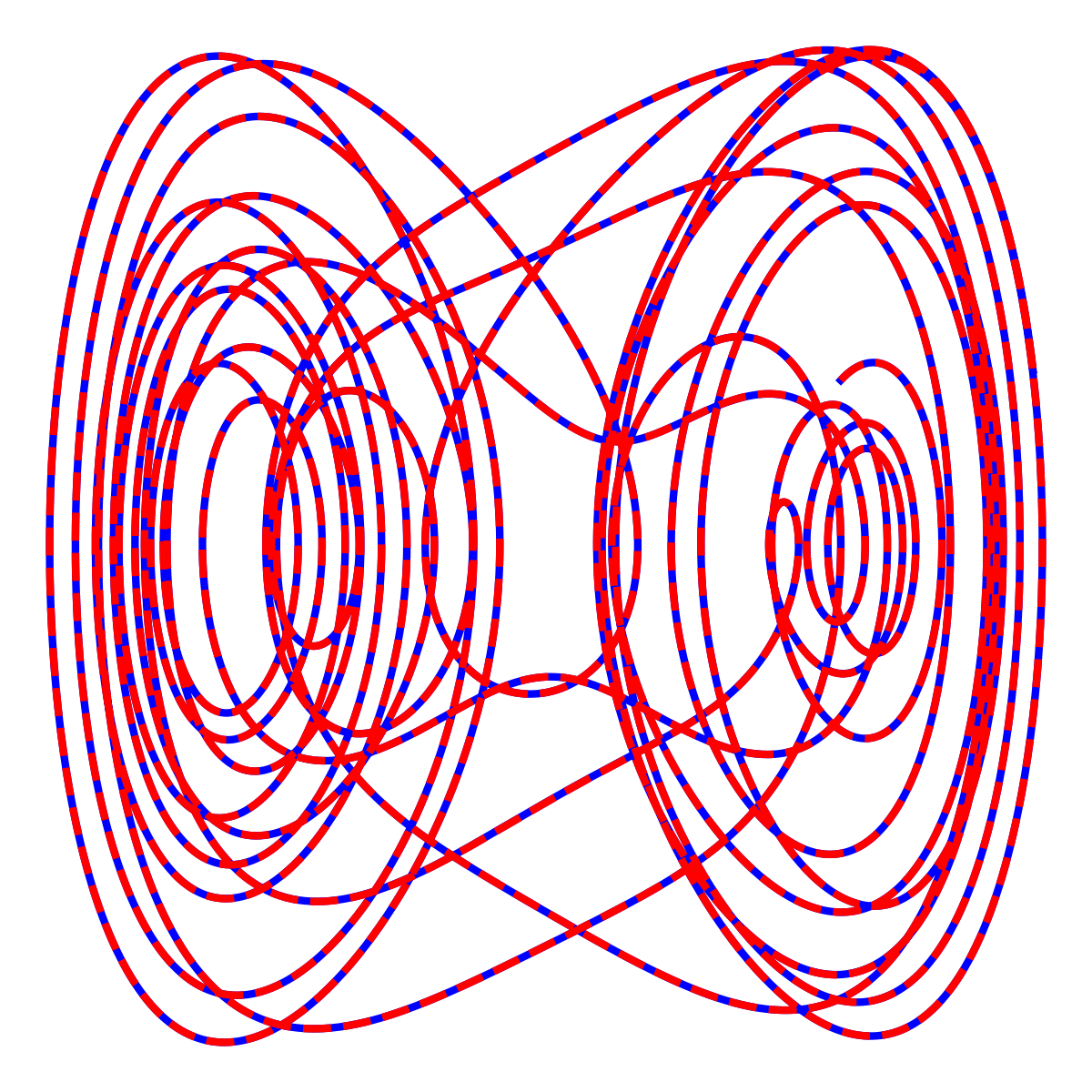}
        \end{overpic}
    };
    \node at ($(dynsys.north)+(0,0.25)$) {\hyperref[sec:dyn_systems]{Dynamical systems}};

    \node at (10,0) (diffop)
    {\begin{overpic}[height=3cm]{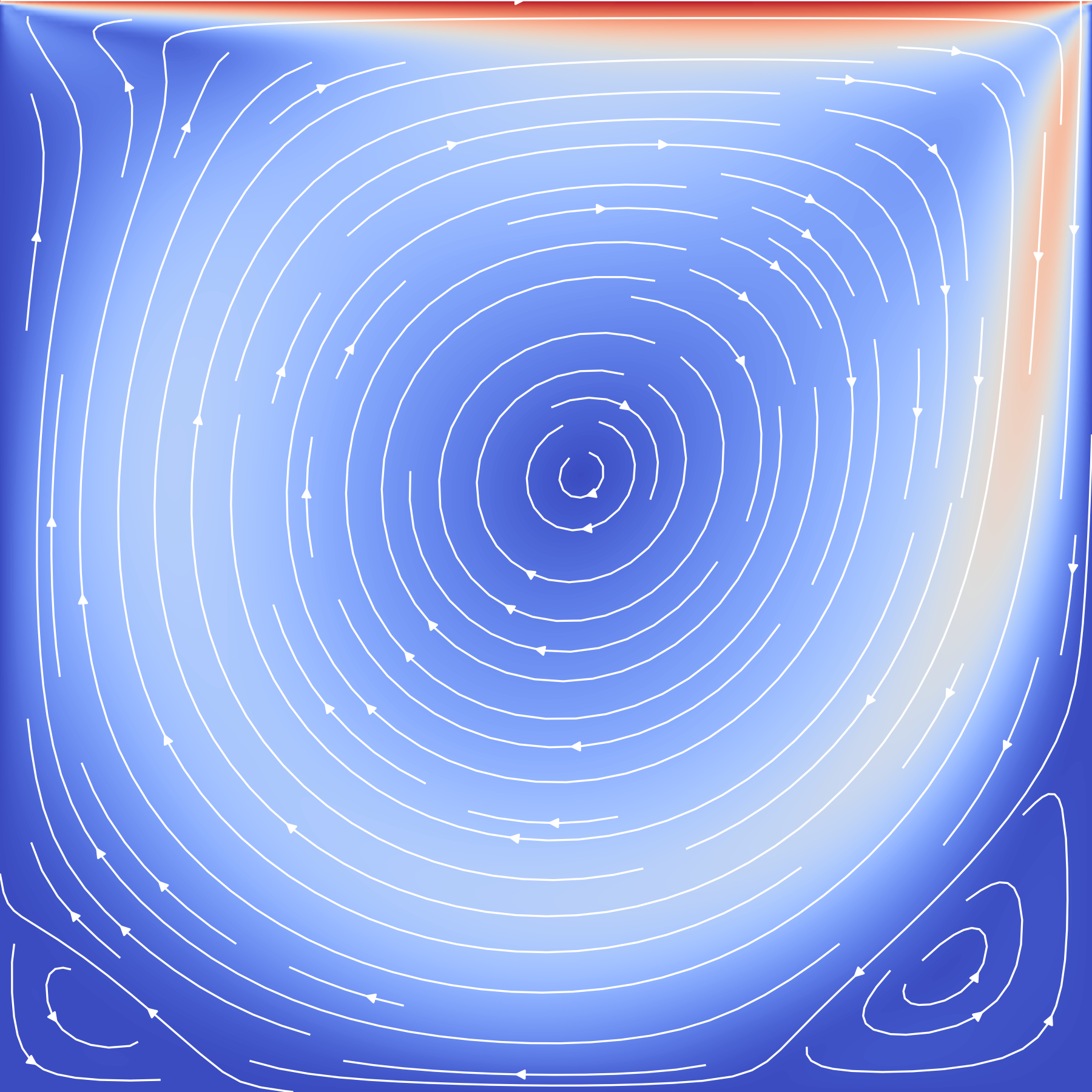}
        \end{overpic}
    };
    \node at ($(diffop.north)+(0,0.25)$) {\hyperref[sec:differential_operators]{Differential operators}};

    \node[box, fill=white] at (5,4) (framework) {\hyperref[sec:analytic_systems]{Analytic systems}\\[5pt]$\sum_{i=1}^N c^*_i \mcal{D}_i(u, f) = \mcal{D}_0(u,f)$};

    \begin{pgfonlayer}{bg}
        \draw[arrow] ($(framework.south west)+(3pt,3pt)$) -- ($(matrix.north)+(0,0.6)$);
    \end{pgfonlayer}
    \draw[arrow] (framework.south) -- ($(dynsys.north)+(0,0.6)$);
    \begin{pgfonlayer}{bg}
        \draw[arrow] ($(framework.south east)+(-3pt,3pt)$) -- ($(diffop.north)+(0,0.6)$);
    \end{pgfonlayer}

\end{tikzpicture}
  }
  \vspace{-0.5cm}
  \caption{We explore three main types of analytic systems in this paper. The core analytic framework (degenerate sets and zero-one law) is developed in \Cref{sec:analytic_systems}, then specialized to matrix recovery (\Cref{sec:matrix_recovery}), analytic dynamical systems (\Cref{sec:dyn_systems}), and differential operators (\Cref{sec:differential_operators}), along with a refined characterization of degenerate input sets for polynomial ODEs and elliptic PDEs in \Cref{sec:refined_characterization}.}
  \label{fig:theory_roadmap}
\end{figure}

Given the importance of model discovery, there are several related lines of work, but they address different identifiability questions. In operator inference, one learns reduced-order models from data generated by known but computationally expensive simulations, typically by solving regression problems in a reduced space rather than by projecting the full governing equations~\cite{peherstorfer2016data,kramer2024learning}. This complements classical projection-based model reduction~\cite{antoulas2005approximation,rozza2008reduced,quarteroni2014reduced,benner2015survey,hesthaven2016certified}. The object  learned in operator inference is usually a reduced surrogate model, whereas our goal is to identify the coefficients of a prescribed differential-operator dictionary. Parameter estimation and structural identifiability are closer in spirit to the present work. In these problems, the form of the governing equation is typically  known, and one asks whether its parameters, or functions of its parameters, are determined by the available observations. Algebraic criteria have been developed for identifying coefficients in ordinary and partial differential equations~\cite{hong2020global,byrne2025algebraic}, and related work characterizes when parameters in linear dynamical systems~\cite{stanhope2014identifiability} or evolution partial differential equations~\cite{he2024much} are recoverable based on a single trajectory.
Random matrix and compressed sensing theory has also been applied to design sampling strategies for recovering sparse solutions to overdetermined linear systems arising from dynamical system identification~\cite{schaeffer2018extracting,schaeffer2020extracting}. Other recent work studies uniqueness of PDEs consistent with a prescribed solution~\cite{scholl2026symbolic}, computes identifiable functions of parameters in ODE models~\cite{ovchinnikov2021computing}, or provides guarantees for discovering chaotic dynamical systems~\cite{shumaylov2025system}. A genericity result related to ours shows for analytic ODEs with $N$-dimensional parameter spaces, a generic $(2N+1)$-tuple of experiments identifies every pair of parameters that can be distinguished by some experiment~\cite{sontag2002differential}. This result applies to nonlinearly-parameterized systems and allows for partial observations, but when it is specialized to an ODE in the form of \cref{eq:model}, it yields an identifiability result using $(2N+1)$ generic experiments.
When the experiment space includes input signals, this does not imply one-shot recovery.

Our focus is complementary. Given a finite dictionary of candidate terms and an unknown model that is some linear combination of those terms, we ask which inputs make the coefficient vector identifiable from one observed input-solution pair, as well as how ``large''  this set of admissible inputs is. Thus, the input itself is central to the identifiability question. For each admissible input, the observed solution induces a finite-dimensional linear map from coefficients to dictionary outputs. We call an input \emph{nondegenerate} when it uniquely identifies the underlying system in one shot, i.e., when this map is injective. Our main result is a zero-one law for analytic systems: if nondegenerate inputs exist, finding them is easy as almost every input sampled from a nondegenerate Gaussian measure is nondegenerate. When one-shot identification is impossible, such an input provides an almost sure certificate of impossibility. For the special case of differential operators in the form $\sum_{i=1}^{N} c_i^* \mcal{D}_i(u) = f$ with dictionary terms depending on the spatial variable and partial derivatives of $u$, we show that nondegenerate inputs always exist, and for univariate polynomial ODEs we find non-smooth deterministic inputs that are always nondegenerate.

The paper is organized as follows. We begin in \Cref{sec:analytic_systems} by introducing the generic analytic framework for studying the prevalence of degenerate input sets for linearly-parameterized operators between Banach spaces. We then apply this framework to matrix recovery problems in \Cref{sec:matrix_recovery}, analytic dynamical systems in \Cref{sec:dyn_systems}, and differential operators in \Cref{sec:differential_operators}. Finally, we provide refined characterizations of the degenerate input sets for certain polynomial ordinary differential equations and partial differential equations in \Cref{sec:refined_characterization}. A schematic overview of the paper is illustrated in \Cref{fig:theory_roadmap}.

A brief announcement of this work was published in \cite{boulle2026zero} without proofs and methods.

\section{Zero-One Law for Analytic Systems Identification}
\label{sec:analytic_systems}
Throughout this section, we consider the system given by
\begin{equation}\label{eqn:analytic_system}
  \mcal{L}(u, f) \coloneqq \sum_{i=1}^N c^*_i \mcal{D}_i(u, f) = \mcal{D}_0(u,f),\qquad (u,f)\in \mcal{U}\times \mcal{F},
\end{equation}
where $\mcal{D}_i : \mcal{U} \times \mcal{F} \to \mcal{V}$ are known analytic functions between Banach spaces \cite{whittlesey1965analytic, bochnak1970analytic, gamelin1994analytic, buffoni2003analytic} and $\vect{c}^*=(c_1^*,\ldots,c_N^*)\in\R^N$ is unknown. We aim to uniquely determine the coefficient vector $\vect{c}^* \in \R^{N}$ given a single input-solution pair $(f,u)$. The formulation in \cref{eqn:analytic_system} allows for spurious dictionary terms: if zero coefficients can be reliably identified for large $N$, then the correct model form can be recovered within dictionaries that include many extra terms. We begin with the following assumption that guarantees that the solution operator $\mcal{G}:\mcal{F}\to\mcal{U}$ associated with the ground truth system \eqref{eqn:analytic_system} is well-posed.

\begin{assumption}[Well-posedness of analytic systems]\label{ass:well_posedness}
  We assume that the solution operator $\mcal{G}: \mcal{F} \to \mcal{U}$ associated with \cref{eqn:analytic_system} is analytic and \emph{well-posed}: for every $f\in \mcal{F}$, there is a unique solution $u \in \mcal{U}$ to \cref{eqn:analytic_system}.
\end{assumption}

Our framework encompasses a wide range of systems, including linearly parameterized families in which the coefficients multiplying the fixed dictionary elements are unknown:
\begin{enumerate}[leftmargin=*, topsep=0pt]
  \item {\bf Boundary value problems.} Here, one seeks to recover a partial differential equation of the form: $\sum_{i=1}^N c^*_i \mcal{D}_i(u) = f$ with dictionary terms depending only on the solution. One example is the Poisson equation on the unit square with homogeneous Dirichlet conditions, i.e.,  $\mcal{L}(u, f) = u_{xx}+u_{yy}$ and $\mcal{D}_0(u, f) = f$ defined from $(H^2(\Omega)\cap H_0^1(\Omega)) \times L^2(\Omega) \to L^2(\Omega)$ on $\Omega=(0,1)^2$. Here, $L^2(\Omega)$ denotes the space of square-integrable functions on $\Omega$, $H^k(\Omega)$ is the Sobolev space of functions with weak derivatives up to order $k$ in $L^2(\Omega)$, and $H_0^1(\Omega)$ denotes the subspace of $H^1(\Omega)$ with zero trace on $\partial\Omega$.
  \item {\bf Dynamical systems.} Consider a dynamical system of the form $\dot{u} = F(t,u,f)$ with initial condition $u(a) = f_a$ for $t\in[a,b]$ (see \Cref{sec:dyn_systems} for more details). One can express this in the form required by~\cref{eq:model} as
        \[\sum_{i = 1}^{N} c^*_i \underbrace{\int_{a}^{t} D_i(s, u(s), f(s)) \d s}_{\mcal{D}_i[u,(f_a, f)](t)} = \underbrace{u(t) - f_a.\vphantom{\int_{a}^{t}}}_{\mcal{D}_0[u,(f_a, f)](t)}\]
        Here, $D_i : [a, b] \times \R^n \times \R^m \to \R^n$ are known, $\mcal{U} = \AC([a,b];\R^n)$  and $\mcal{F} \subset \R^n \times L^{\infty}([a,b]; \R^m)$, where $L^\infty([a,b];\mathbb{R}^m)$ denotes the space of essentially bounded functions from $[a,b]$ to $\mathbb{R}^m$, and $\AC([a,b];\mathbb{R}^n)$ denotes the space of absolutely continuous functions from $[a,b]$ to $\mathbb{R}^n$.
  \item {\bf Matrix recovery.} Our framework also applies to families of models that are not differential operators. To demonstrate this, we consider matrix recovery problems where one aims to recover a linearly-parameterized matrix of the form $A = \sum_{i=1}^{N} c^*_i A_i$, where $A_i\in\R^{m\times n}$, from matrix-vector products.
\end{enumerate}

These examples are not exhaustive. The only structural requirements in the abstract theory are analytic dependence on the observed input-solution pair and finite-dimensional linear dependence on the unknown coefficients.

\subsection{Genericity of one-shot identification} \label{sec:genericity}

One-shot identification of~\eqref{eqn:analytic_system} is not always possible. If there exists an input-output pair $(f,u)$ such that the dictionary elements $\mathcal{D}_i(u,f)\in \mathcal{V}$ are linearly dependent, i.e., there is a nonzero vector of coefficients $\vect{a}\in\R^N$ such that $\sum_{i=1}^{N} a_i \mathcal{D}_i(u,f)=0$, then there are at least two distinct choices ($\vect{c}^*$ and $\vect{c}^*+\vect{a}$) that exactly describe the input-solution pair $(u,f)$. If we are unfortunate enough to observe that particular input-solution pair, then we  cannot uniquely identify the coefficients in the model, and at least one more input-solution pair would be required.

Therefore, to understand one-shot identification, we begin by understanding the injectivity of the so-called parameter-to-output recovery map $L_{f, u}:\R^N \to \mcal{V}$ associated with~\cref{eqn:analytic_system}, defined for a fixed $(f, u)\in \mathcal{F}\times \mathcal{U}$ by the following linear map:
\begin{equation}\label{eqn:analytic_system_recoverability_operator}
  L_{f, u}: \vect{c} \longmapsto \sum_{i=1}^N c_i \mcal{D}_i(u, f).
\end{equation}
If $L_{f,u}$ is injective, then there is a unique vector of coefficients such that $L_{f, u} [\vect{c}] = \mcal{D}_0(u, f)$. Existence of such a vector is guaranteed for data generated by~\cref{eqn:analytic_system}. On the other hand, when $L_{f, u}$ is not injective, one-shot identification is not possible using that particular $(f,u)$ pair.

Thus, to understand when one-shot identification is not possible, we seek to identify the set of inputs $f\in\mathcal{F}$ for which the operator $L_{f, \mcal{G}(f)}$ in~\cref{eqn:analytic_system_recoverability_operator} is not injective, where $\mcal{G}:\mathcal{F}\to\mathcal{U}$ is the solution operator associated with \cref{eqn:analytic_system}. This leads us to the notion of degenerate input sets.

\begin{definition}[Degenerate input set]\label{def:degenerate_set}
  We define the degenerate set of inputs associated with~\cref{eqn:analytic_system} as
  \[
    \Fdeg(\vect{c}^*) =
    \left\{ f \in \mcal{F} \mid \exists \vect{a}\in\R^N\setminus\{0\},\, L_{f, \mcal{G}(f)}(\vect{a}) = 0  \right\},
  \]
  where $\mcal{G}$ is the solution operator for the ground-truth coefficients $\vect{c}^*$.
\end{definition}
The set $\Fdeg(\vect{c}^*)\subset \mcal{F}$ characterizes the input functions that make the operator $L_{f, \mcal{G}(f)}$ in~\cref{eqn:analytic_system_recoverability_operator} not injective, and these are the inputs that we would like to avoid, as our goal is one-shot identification. To avoid the set $\Fdeg(\vect{c}^*)$, we need to understand its ``prevalence" in $\mcal{F}$. First, \cref{def:compactly_slice_null} introduces the notion of compactly slice-null sets, which characterizes sets that have zero Lebesgue measure when intersected with sufficiently large finite-dimensional affine subspaces.

\begin{figure}[t]
  \centering
  \begin{tikzpicture}[x=0.03636\textwidth,y=0.03636\textwidth]
    \begin{scope}
        %\draw[style=help lines,step=1cm] (-4,-4) grid (4,4);

        \draw[->,thick] (-4.5,0) -- (4.5,0);
        \node[anchor=north] at (4.5,0) {$x$};
        \draw[->,thick] (0,-3) -- (0,3);
        \node[anchor=west] at (0,3) {$y$};

        \draw[very thick, red, name path=curve 1] (-4.5,-2) -- (4.5,-2);
        \node[black, above left] at (-2,2) {$\phi^{-1}(0)$};
        \node[black, below right] at (2,-2) {$f_1 + E_n$};
        \draw[very thick, blue, name path=curve 2]
        plot [smooth, tension=0.9] coordinates {
                (-4,1.2)
                (-2,2.2)
                (0,2.4)
                %(1,2.0)
                (1.5,1.5)
                (1.0,0.0)
                (0.2,-1.8)
                (-1.0,-2.2)
                (-2.0,-1.0)
                (-1.5,0.5)
                %(-0.5,1.0)
                (1.0,-0.8)
                (2.0,0.6)
                (4,0.4)
            };

        \tikzset{name intersections={of=curve 1 and curve 2, name=i, total=\t}}

        \fill[green!70!black, every node/.style={above left, black, opacity=1}]
        \foreach \s in {1,...,\t}{(i-\s) circle (3pt) node {}};
    \end{scope}

    % Center cue indicating vertical translates of E_n (midpoint between panels)
    \draw[->, very thick, red] (6.7,-2.7) -- (6.7,2.7);
    \node[rotate=90, black] at (6,0) {translates of $E_n$};

    \begin{scope}[shift={(13,0)}]
        \draw[->,thick] (-4.5,0) -- (4.5,0);
        \node[anchor=north] at (4.5,0) {$x$};
        \draw[->,thick] (0,-3) -- (0,3);
        \node[anchor=west] at (0,3) {$y$};

        \draw[very thick, red, name path=curve 1] (-4.5,1.5) -- (4.5,1.5);
        \node[black, above right] at (2,1.5) {$f_2+ E_n$};
        \node[black, above left] at (-2,2) {$\phi^{-1}(0)$};
        \draw[very thick, blue, name path=curve 2]
        plot [smooth, tension=0.9] coordinates {
                (-4,1.2)
                (-2,2.2)
                (0,2.4)
                %(1,2.0)
                (1.5,1.5)
                (1.0,0.0)
                (0.2,-1.8)
                (-1.0,-2.2)
                (-2.0,-1.0)
                (-1.5,0.5)
                %(-0.5,1.0)
                (1.0,-0.8)
                (2.0,0.6)
                (4,0.4)
            };

        \tikzset{name intersections={of=curve 1 and curve 2, name=i, total=\t}}

        \fill[green!70!black, every node/.style={above left, black, opacity=1}]
        \foreach \s in {1,...,\t}{(i-\s) circle (3pt) node {}};
    \end{scope}
\end{tikzpicture}
  \caption{Illustration of \cref{thm:analytic_varieties_are_small}. We show that the zero level set $\phi^{-1}(0)$ (blue curve) of a nonzero analytic function $\phi:\mcal{F} \to \R$ is
    \emph{compactly slice-null} in the sense that for any compact set, there is a sufficiently large finite-dimensional probe space $E_n$ such that all translates of the probe space by points $f$ in the compact set (red line) have Lebesgue null intersection with $\phi^{-1}(0)$ (green dots). By a Fubini-like theorem, the analytic variety $\phi^{-1}(0)$ has zero Gaussian measure.}
  \label{fig:slices_of_analytic_variety}
\end{figure}
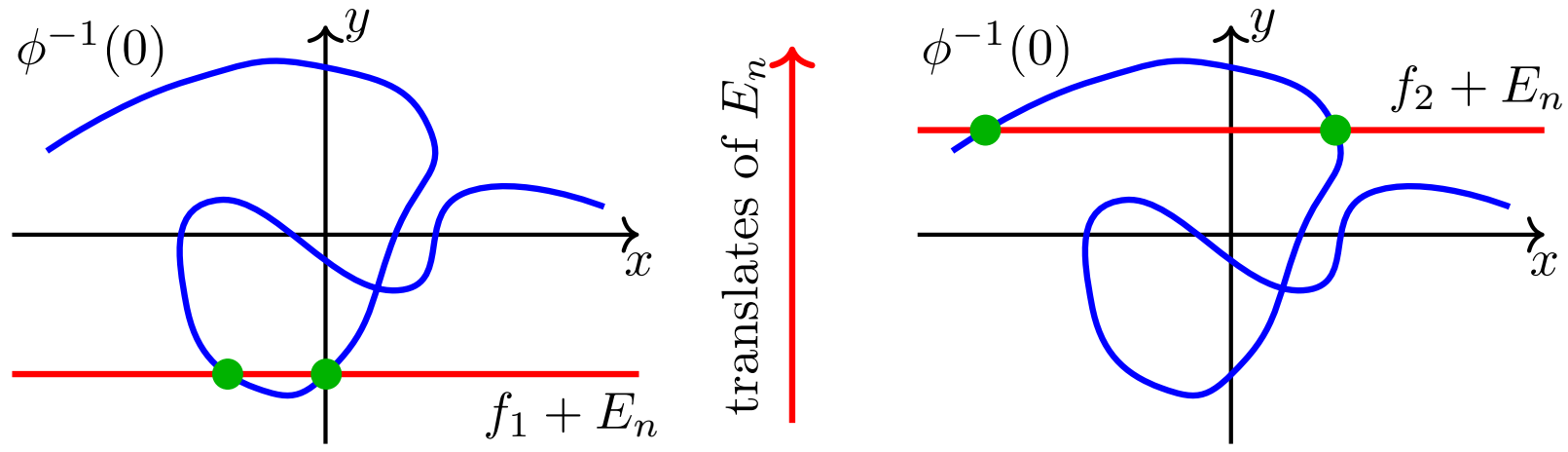

\begin{definition}[Compactly slice-null set]\label{def:compactly_slice_null}
  Let $\mcal{F}$ be a Banach space with a countable Schauder basis. A set $S \subset \mcal{F}$ is \emph{compactly slice-null} if for every countable Schauder basis $\mcal{E} = \{e_j\}_{j=1}^{\infty}$ and every compact set $K \subset \mcal{F}$, there exists $N_{K, \mcal{E}} \geq 1$ such that for every $n \geq N_{K, \mcal{E}}$ and every $f \in K$, the affine slice $S \cap (f + E_n)$ has zero Lebesgue measure in the finite-dimensional affine subspace $f + E_n$, where $E_n = \vspan\{ e_j \}_{j=1}^{n}$.
\end{definition}

We now prove \cref{thm:analytic_varieties_are_small} on analytic varieties in Banach spaces, which is illustrated in \cref{fig:slices_of_analytic_variety}. A Gaussian measure $\gamma$ on $\mcal{F}$ is called nondegenerate if every nonzero continuous linear functional on $\mcal{F}$ has positive variance under $\gamma$. This generalizes a standard result that the zero level set of a nonzero analytic function in finite dimensions has zero Lebesgue measure \cite[p.~240]{federer1969geometric} to functions defined on infinite-dimensional Banach spaces. The first half of the theorem provides a fine characterization of the zero level set when restricted to finite-dimensional slices and allows us to recover a previous result by Bogachev \& Malofeev~\cite[Cor.~2]{bogachev2014distributions} in \cref{thm:analytic_Gaussian_null}.

\begin{theorem}[Zero-level sets of analytic functions]\label{thm:analytic_varieties_are_small}
  Let $\mcal{F}$ be a Banach space with a countable Schauder basis. If $\phi: \mcal{F} \to \R$ is a nonzero analytic function, then the following two results hold:
  \begin{enumerate}[label=\roman*), itemsep=0em, leftmargin=*]
    \item \label{thm:analytic_shyness} The zero level set of $\phi$ is compactly slice-null.
    \item \label{thm:analytic_Gaussian_null}$\phi^{-1}(0)$ has measure zero with respect to any nondegenerate Gaussian measure on Borel sets of $\mcal{F}$.
  \end{enumerate}
\end{theorem}
\begin{proof}
  We begin by proving \cref{thm:analytic_shyness}, and then we use it  to prove \cref{thm:analytic_Gaussian_null}. Let $\mcal{E} = \{e_j\}_{j=1}^{\infty}$ be a countable Schauder basis for $\mcal{F}$ and let $E_n = \vspan\{ e_j \}_{j=1}^{n}$ for $n\geq 1$. Let $f_0 \in \mcal{F}$ be such that $\phi(f_0) \neq 0$, and consider an arbitrary $f' \in K$.
  Since $\mcal{E}$ is a Schauder basis, we can find a sequence $\{\Delta f_n\}_{n=1}^{\infty}$ with $\Delta f_n \in E_n$ and $\Delta f_n \to f_0 - f'$ as $n\to\infty$.
  Hence, $f' + \Delta f_n \in f' + E_n$ and $f' + \Delta f_n \to f_0$ as $n \to \infty$.
  Since $\phi(f_0) \neq 0$ and $\phi$ is continuous, there exists $N_{f'} \geq 1$ such that $\Delta \tilde{f} \coloneqq \Delta f_{N_{f'}} \in E_{N_{f'}}$ satisfies  $\phi(f' + \Delta \tilde{f}) \neq 0$.
  Furthermore, by continuity of $\phi$, there exists an $R_{f'} > 0$ such that for every $f \in \mcal{F}$ with $\| f - f' \| < R_{f'}$ we have $\phi(f + \Delta \tilde{f}) \neq 0$.
  Since the balls  $B(f', R_{f'}) = \{ f \in \mcal{F} \ : \ \| f - f' \| < R_{f'} \}$ with $f' \in K$ form an open cover of the compact set $K$, there is a finite subcover  $K \subset \bigcup_{i=1}^{M} B(f_i, R_{f_i})$, and corresponding $N_{f_i}$ with $\Delta \tilde{f}_i \in E_{N_{f_i}}$ satisfying $\phi(f + \Delta \tilde{f}_i) \neq 0$ for all $f \in B(f_i, R_{f_i})$.

  Setting $N_K = \max\{ N_{f_1}, \ldots, N_{f_M} \}$, choose any $n \geq N_K$ and $f \in K$.
  Since $\phi$ is analytic, the mapping $\phi_{f,n}: E_n \to \R$ defined by
  $\phi_{f,n}: \Delta f \mapsto \phi(f + \Delta f)$
  is analytic and satisfies $
    \phi^{-1}(0) \cap (f + E_n)
    = f + \phi_{f,n}^{-1}(0)$.
  It remains to show that $f + \phi_{f,n}^{-1}(0)$ has zero Lebesgue measure in $f + E_n$. Since $f \in K$, there is an index $i\in\{1,\ldots,M\}$ such that $f \in B(f_i, R_{f_i})$. Moreover, since $n \geq N_K \geq N_{f_i}$, we have $\Delta \tilde{f}_i \in E_{N_{f_i}} \subset E_n$ and
  $\phi_{f,n}(\Delta \tilde{f}_i) = \phi(f + \Delta \tilde{f}_i) \neq 0$. If $\phi_{f,n}$ is constant on $E_n$, then $\phi_{f,n}^{-1}(0)$ is empty.
  Otherwise, if $\phi_{f,n}$ is not constant on $E_n$, then it follows from
  the result stated on \cite[p.~240]{federer1969geometric}, that $\phi_{f,n}^{-1}(0)$ has zero Lebesgue measure in $E_n$.
  Therefore, $f + \phi_{f,n}^{-1}(0)$ has zero Lebesgue measure in $f + E_n$, completing proof of \cref{thm:analytic_shyness}.

  Turning to the proof of \cref{thm:analytic_Gaussian_null},
  let $\gamma$ be a nondegenerate Gaussian measure on Borel sets of $\mcal{F}$, noting that the cylindrical, Baire, and Borel $\sigma$-algebras all coincide in the separable Banach space $\mcal{F}$ by \cite[Chapt.~I, Thm.~1.2]{vakhania1987probability}.
  It follows from \cite[Chapt.~I, Thm.~3.1]{vakhania1987probability} that $\gamma$ is a Radon measure, hence it is tight. Let $\varepsilon > 0$. By tightness of $\gamma$, there exists a compact subset $K_{\varepsilon} \subset \mcal{F}$ such that $\gamma(\mcal{F} \setminus K_{\varepsilon}) < \varepsilon$.
  By \cite[Thm.~3.6.1]{bogachev1998gaussian}, nondegeneracy of $\gamma$ implies that its Cameron--Martin space $H(\gamma)$ is dense in $\mcal{F}$.
  Due to the stability of Schauder bases given by \cite[Prop.~1.a.9(i)]{lindenstrauss1996classical}, it follows that $\mcal{F}$ admits a (perturbed) Schauder basis $\{ \tilde{e}_n \}_{n=1}^{\infty} \subset H(\gamma)$.
  Since $\phi$ is continuous, $\phi^{-1}(0)$ is Borel measurable, and it suffices to show that $\gamma(\phi^{-1}(0)) = 0$.
  By \cref{thm:analytic_shyness}, there is an integer $N_{\varepsilon}$ such that with $E = \vspan\{ \tilde{e}_j \}_{j=1}^{N_{\varepsilon}} \subset H(\gamma)$, each slice $\phi^{-1}(0) \cap(f + E)$ for $f\in K_{\varepsilon}$ has zero Lebesgue measure in the finite-dimensional affine subspace $f + E$.

  Since $E$ is finite-dimensional, it admits a closed complementary space $Y$ such that $\mcal{F} = E \oplus Y$ with bounded linear projections onto each component \cite[Sec.~2.4]{brezis2010functional}.
  With $\Pi$ being the projection onto $Y$, \cite[Cor.~3.10.3]{bogachev1998gaussian} yields a measure $\gamma^y$ on each fiber $y + E = \Pi^{-1}(y)$, $y\in Y$, such that
  \[
    \gamma\big( \phi^{-1}(0) \cap K_{\varepsilon} \big)
    = \int_{Y} \gamma^{y}\big( \phi^{-1}(0) \cap K_{\varepsilon} \big) \Pi_*\gamma(\rm{d} y).
  \]
  For $\Pi_*\gamma$-almost every $y\in Y$, the conditional measure $\gamma^y$ is equivalent to Lebesgue measure on the finite-dimensional affine space $y+E$. For such a $y$, if $(y+E)\cap K_\varepsilon=\emptyset$, then $\gamma^y(\phi^{-1}(0)\cap K_\varepsilon)=0$. Otherwise, choose $f\in(y+E)\cap K_\varepsilon$. Then $y+E=f+E$, and
  \[
    \phi^{-1}(0) \cap K_{\varepsilon} \cap (y + E) \subset \phi^{-1}(0) \cap(f + E),
  \]
  which has zero Lebesgue measure in $y+E$ by the slice-null conclusion. Hence $\gamma^{y}( \phi^{-1}(0) \cap K_{\varepsilon} ) = 0$ for $\Pi_*\gamma$-almost every $y$, and therefore $\gamma ( \phi^{-1}(0) \cap K_{\varepsilon} ) = 0$.
  Since
  \[
    \gamma\big( \phi^{-1}(0) \big)
    = \gamma\big( \phi^{-1}(0) \cap K_{\varepsilon} \big) + \gamma\big( \phi^{-1}(0) \setminus K_{\varepsilon} \big)
    < \varepsilon
  \]
  for every $\varepsilon > 0$, it follows that $\gamma( \phi^{-1}(0) ) = 0$, completing the proof of \cref{thm:analytic_Gaussian_null}.
\end{proof}

In the following counterexample, we show that it is essential to use a non-degenerate Gaussian measure  in \cref{thm:analytic_varieties_are_small}.

\begin{example}[Degenerate Gaussian measure]
  Consider the function $\phi:\R^2\to\R$ defined as $\phi(x,y) = y$ for all $(x,y)\in \R^2$. Let $\mu$ be a Gaussian measure on $\R$, $\delta_0$ the Dirac measure at $y=0$, and $\gamma = \mu\otimes \delta_0$ be a degenerate Gaussian measure with rank-$1$ covariance kernel. Then, $\phi^{-1}(0) = \vspan\{(1,0)\} \neq \R^2$, but $\gamma(\phi^{-1}(\{0\})) = \mu(\R) = 1$.
\end{example}

The next theorem shows that if the coefficient vector $\vect{c}^*$ can be recovered from an input $f_0\in\mathcal{F}$ and its solution $u_0 = \mcal{G}(f_0)$, then it can also be recovered using ``generic'' choices of $f$, sampled from a Gaussian measure on $\mathcal{F}$, and corresponding solutions $u = \mcal{G}(f)$.

\begin{theorem}[Prevalence of degenerate sets] \label{thm:degenerate_sets}
  Let $\mcal{V}$ and $\mcal{F}$ be Banach spaces such that $\mcal{F}$ has a countable Schauder basis, and let $\mcal{B}(\R^N,\mcal{V})$ denote the Banach space of bounded linear operators $\R^N \to \mcal{V}$. Consider a continuous mapping $f\in\mcal{F} \mapsto L_f\in \mcal{B}(\R^N,\mcal{V})$ such that for every $\varphi \in \mcal{V}^*$ and $\vect{v} \in \R^N$ the real-valued mapping $f \mapsto \varphi(L_f \vect{v})$ is analytic, with associated degenerate set $\Fdeg = \left\{ f \in \mcal{F} \mid \rank(L_f) \leq N-1 \right\}$. If there exists $f_0 \in \mcal{F} \setminus \Fdeg$, then $\Fdeg$ is compactly slice-null in $\mcal{F}$ and has measure zero with respect to any nondegenerate Gaussian measure on Borel sets of $\mcal{F}$.
\end{theorem}
\begin{proof}
  Let $1\leq k\leq N$. We will prove a stronger result for the set $\Fdeg^k = \left\{ f \in \mcal{F} \mid \rank(L_f) \leq N-k \right\}$, noting that $\Fdeg^1 = \Fdeg$. First, we establish that $\Fdeg^k$ is measurable with respect to the cylindrical, Baire, and Borel $\sigma$-algebras, which all coincide in the separable Banach space $\mcal{F}$ \cite[Chapt.~I, Thm.~1.2]{vakhania1987probability}. Let $S$ be a countable dense set in the Stiefel manifold $\St(N,k)$ consisting of real $N \times k$ matrices with orthonormal columns, and let $\psi: \mcal{F} \to \R$ be defined as $
    \psi(f)
    = \inf_{\mat{V} \in S} \| L_{f} \mat{V} \|$ for $f\in\mathcal{F}$.
  As a countable infimum of continuous functions, $\psi$ is upper semicontinuous, and thus Borel measurable.
  Moreover, we have $\Fdeg^k = \psi^{-1}(0)$, which proves this set is measurable.
  To establish that $\Fdeg^k = \psi^{-1}(0)$, we first observe that by the rank-nullity theorem, if $f\in \Fdeg^k$ there are $k$ linearly independent coefficient vectors $\vect{c}_1, \ldots, \vect{c}_k \in \R^N$ satisfying $L_{f}\vect{c}_j = 0$, which implies that $L_{f}\mat{\tilde{V}} = 0$, where the columns of $\mat{\tilde{V}} \in \St(N,k)$ have the same span as $\vect{c}_1, \ldots, \vect{c}_k$.
  By density of $S$ in $\St(N,k)$, there is a sequence $\{\mat{V}_n\}_{n=1}^{\infty} \subset S$ with $\mat{V}_n\to\mat{\tilde{V}}$.
  By continuity of $L_{f}$, we have
  \[
    \inf_{\mat{V}\in S}\|L_{f}\mat{V}\|
    \leq \|L_{f}\mat{V}_n\| \xrightarrow[n\to\infty]{} \|L_{f}\mat{\tilde{V}}\| = 0,
  \]
  which implies that $\Fdeg^k \subset \psi^{-1}(0)$.
  For the other inclusion, we let $f\in \psi^{-1}(0)$ such that $\inf_{\mat{V} \in S} \| L_{f} \mat{V} \| = 0$.
  Then, there exists a sequence $\{\mat{V}_n\}_n\subset S$ such that $\| L_{f} \mat{V}_n \| \xrightarrow[n\to\infty]{} 0$.
  Since the Stiefel manifold $\St(N,k)$ is compact, there exists a convergent subsequence $\mat{V}_{n_k}\xrightarrow[k\to\infty]{} \mat{V}\in \St(N,k)$, such that, by boundedness of $L_{f}$ and continuity of the operator norm, $\|L_{f} \mat{V}\| = \lim_{k\to\infty}\|L_{f} \mat{V}_{n_k}\| = 0$.
  Since this means $\Range(\mat{V}) \subset \Null L_{f}$, the rank and nullity theorem yields $\rank(L_{f}) \leq N - k$ and so $f \in \Fdeg^k$.
  Therefore, $\psi^{-1}(0)\subset \Fdeg^k$.

  By \cref{thm:analytic_varieties_are_small}, it suffices to construct a nonzero analytic function $\phi : \mcal{F} \to \R$ such that $\Fdeg^k \subset \phi^{-1}(0)$.
  Since $f_0 \notin \Fdeg^k$, we observe that $W = \Range(L_{f_0})$ is a finite-dimensional subspace of $\mcal{V}$ with dimension $r_0 \coloneqq \dim W > N-k$.
  Let $\tilde{\psi}_1, \ldots, \tilde{\psi}_{r_0}$ be a basis for the dual space of $W$, which we denote $W^*$.
  By the Hahn--Banach theorem, these can be extended to elements $\psi_1, \ldots, \psi_{r_0} \in \mcal{V}^*$ which define a bounded linear operator $\Psi: \mcal{V} \to \R^{r_0}$ given by
  \[
    \Psi: g \mapsto \left( \psi_1(g), \ldots, \psi_{r_0}(g) \right), \quad g \in \mcal{V}.
  \]
  Since $L_{f_0}$ has rank $r_0$, we can find columns $e_{j_1}, \ldots, e_{j_{r_0}}$ of the $N\times N$ identity matrix such that $L_{f_0}e_{j_1}, \ldots, L_{f_0} e_{j_{r_0}}$ are linearly independent.
  Let $E \in \R^{N \times r_0}$ be the submatrix of the $N \times N$ identity matrix formed by these columns.
  $L_{f_0} E$ is a linear isomorphism onto $W$ and $\Psi|_W$ is a linear isomorphism from $W\to \R^{r_0}$.
  Hence, the composition $\Psi L_{f_0} E$ is also a linear isomorphism from $\R^{r_0}\to\R^{r_0}$ that can be represented by an invertible real $r_0 \times r_0$ matrix.
  We define $\phi: \mcal{F} \to \R$ by
  \begin{equation} \label{eqn:phi}
    \phi(f)
    = \det\left( \Psi L_{f} E \right)
    = \det\begin{bmatrix}
      \psi_1 (L_{f} e_{j_1})     & \cdots & \psi_1 (L_{f} e_{j_{r_0}})     \\
      \vdots                     & \ddots & \vdots                         \\
      \psi_{r_0} (L_{f} e_{j_1}) & \cdots & \psi_{r_0} (L_{f} e_{j_{r_0}})
    \end{bmatrix}.
  \end{equation}
  As a composition of analytic functions, $\phi$ is analytic \cite{whittlesey1965analytic} and nonzero since $\phi(f_0) \neq 0$ because the matrix $\Psi L_{f_0} E$ is invertible.
  To show that $\phi(\Fdeg^k) = \{ 0 \}$, suppose $f \in \mcal{F}$ with $\phi(f) \neq 0$.
  Then, we have $N - k < r_0 = \rank\left(\Psi L_{f} E \right)
    \leq \rank \left(L_{f}\right)$,
  meaning $f \notin \Fdeg^k$.
  Therefore, if $f \in \Fdeg^k$, then $\phi(f) = 0$, i.e., $\Fdeg^k \subset \phi^{-1}(0)$.
  Applying \cref{thm:analytic_varieties_are_small} completes the proof.
\end{proof}

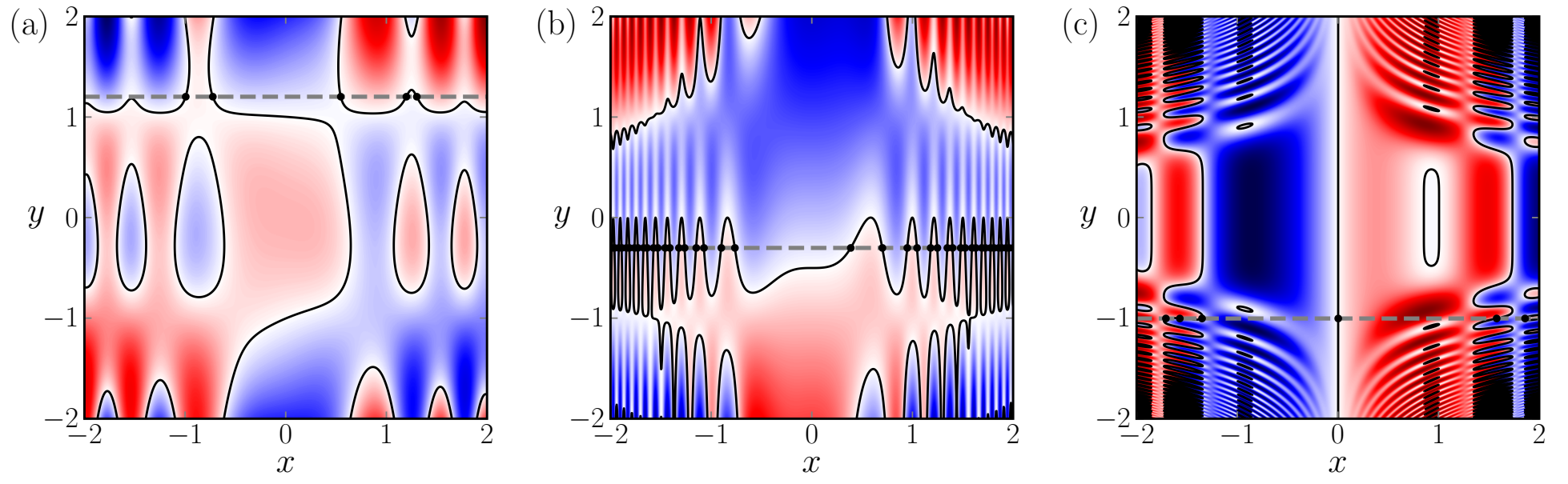
\begin{figure}[t]
  \centering
  \resizebox{\textwidth}{!}{
    \newcommand{\tickfontsize}{\LARGE}
\newcommand{\labelfontsize}{\huge}

\begin{tikzpicture}
    \begin{axis}[
            name=axisA,
            at={(-9cm,1.5cm)},
            anchor=south west,
            width=10cm,
            height=10cm,
            xmin=-2, xmax=2,
            ymin=-2, ymax=2,
            xtick={-2,-1,0,1,2},
            ytick={-2,-1,0,1,2},
            %axis lines=none,
            %ticks=none,
            clip=true,
            axis on top,
            axis line style={line width=1.4pt},
            tick style={line width=0.9pt},
            tick label style={font=\tickfontsize},
            xlabel={$x$},
            ylabel={$y$},
            label style={font=\labelfontsize},
            ylabel style={rotate=-90, xshift=0.3cm},
            axis equal image,
            %scale only axis
        ]
        \addplot graphics [xmin=-2, xmax=2, ymin=-2, ymax=2] {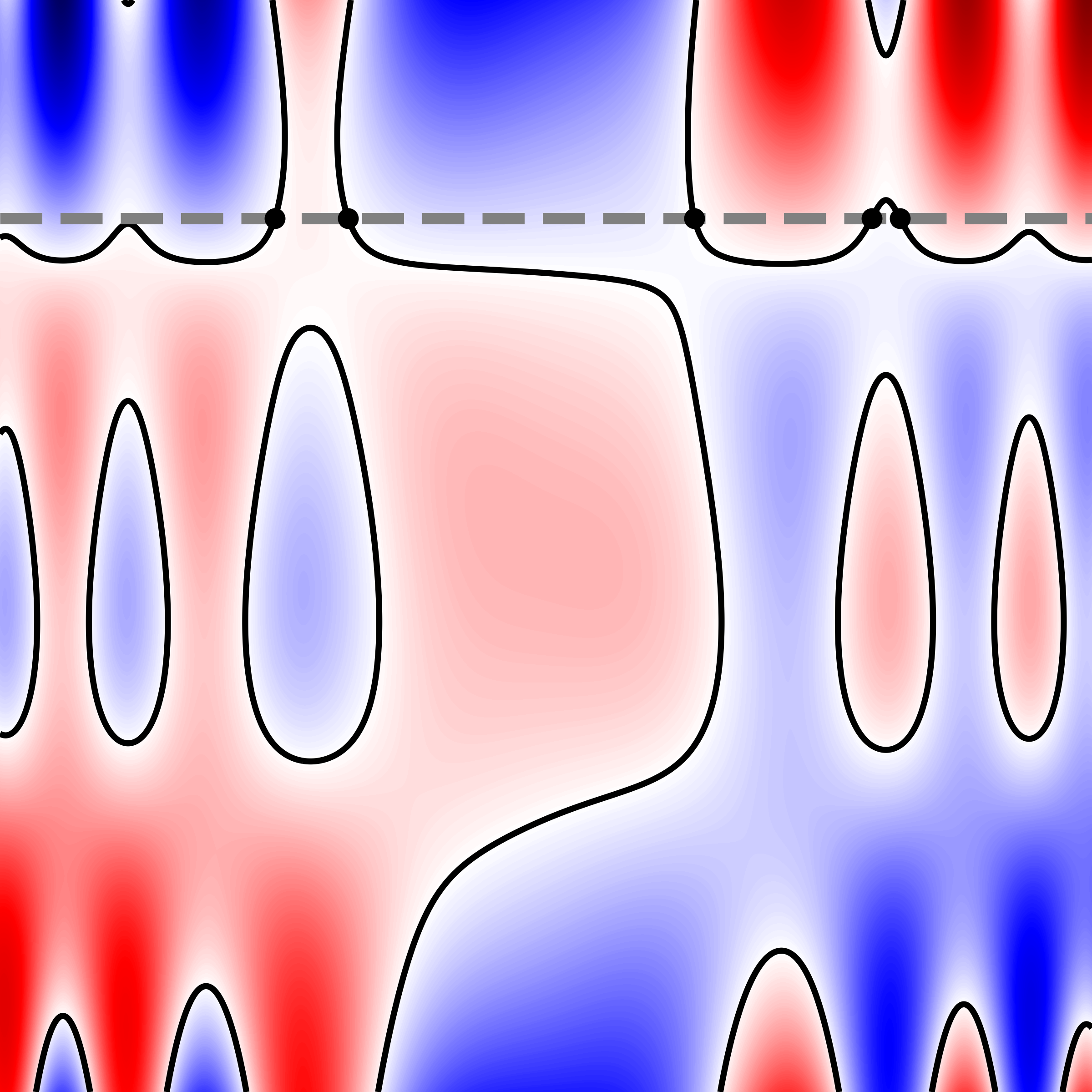};

    \end{axis}

    \begin{axis}[
            name=axisB,
            at={(-9cm,1.5cm)},
            anchor=south west,
            width=10cm,
            height=10cm,
            xmin=-2, xmax=2,
            ymin=-2, ymax=2,
            xtick={-2,-1,0,1,2},
            ytick={-2,-1,0,1,2},
            %axis lines=none,
            %ticks=none,
            clip=true,
            axis on top,
            axis line style={line width=1.4pt},
            tick style={line width=0.9pt},
            tick label style={font=\tickfontsize},
            xlabel={$x$},
            ylabel={$y$},
            label style={font=\labelfontsize},
            ylabel style={rotate=-90, xshift=0.3cm},
            axis equal image,
            xshift=11cm,
            %scale only axis
        ]
        \addplot graphics [xmin=-2, xmax=2, ymin=-2, ymax=2] {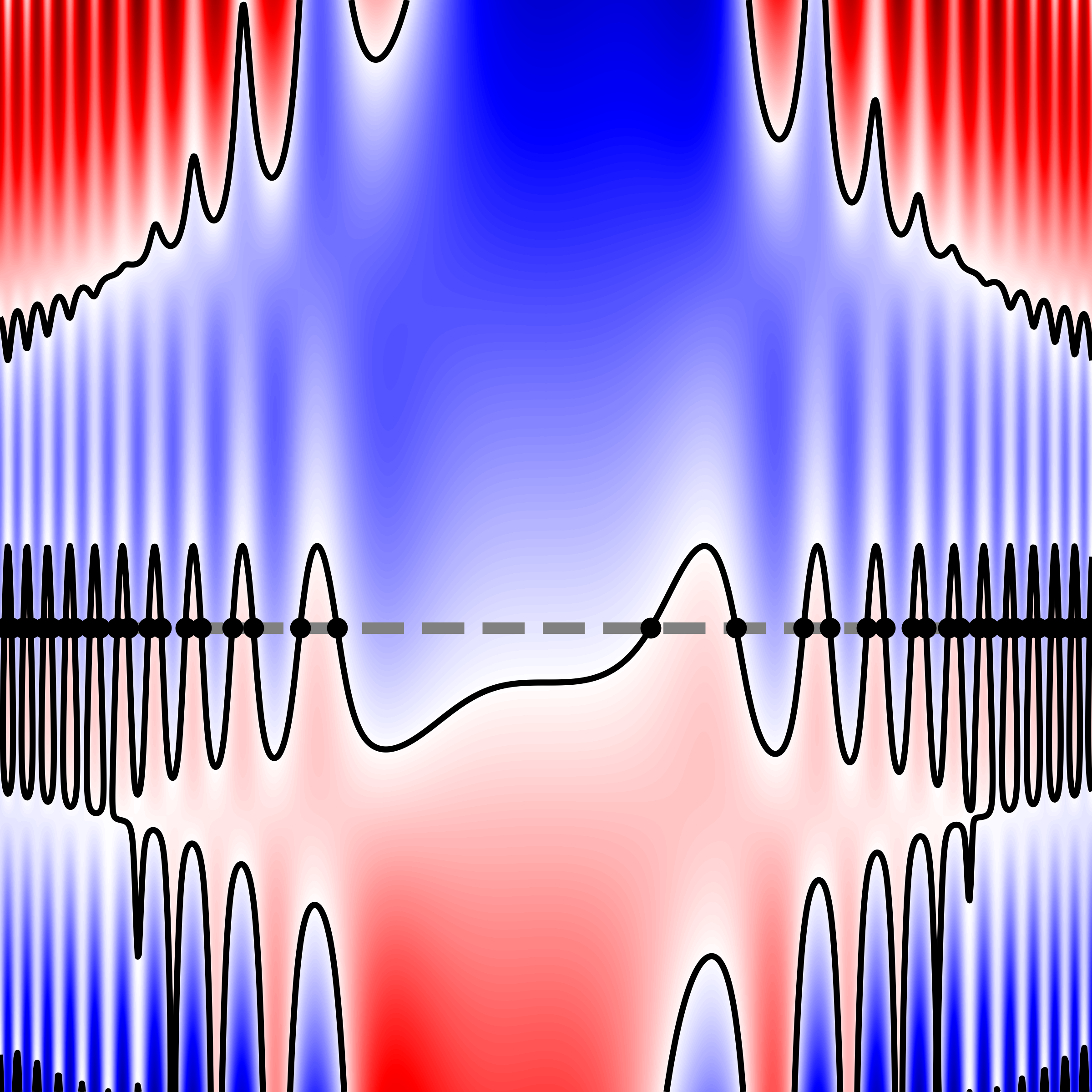};

    \end{axis}

    \begin{axis}[
            name=axisC,
            at={(-9cm,1.5cm)},
            anchor=south west,
            width=10cm,
            height=10cm,
            xmin=-2, xmax=2,
            ymin=-2, ymax=2,
            xtick={-2,-1,0,1,2},
            ytick={-2,-1,0,1,2},
            %axis lines=none,
            %ticks=none,
            clip=true,
            axis on top,
            axis line style={line width=1.4pt},
            tick style={line width=0.9pt},
            tick label style={font=\tickfontsize},
            xlabel={$x$},
            ylabel={$y$},
            label style={font=\labelfontsize},
            ylabel style={rotate=-90, xshift=0.3cm},
            axis equal image,
            xshift=22cm,
            %scale only axis
        ]
        \addplot graphics [xmin=-2, xmax=2, ymin=-2, ymax=2] {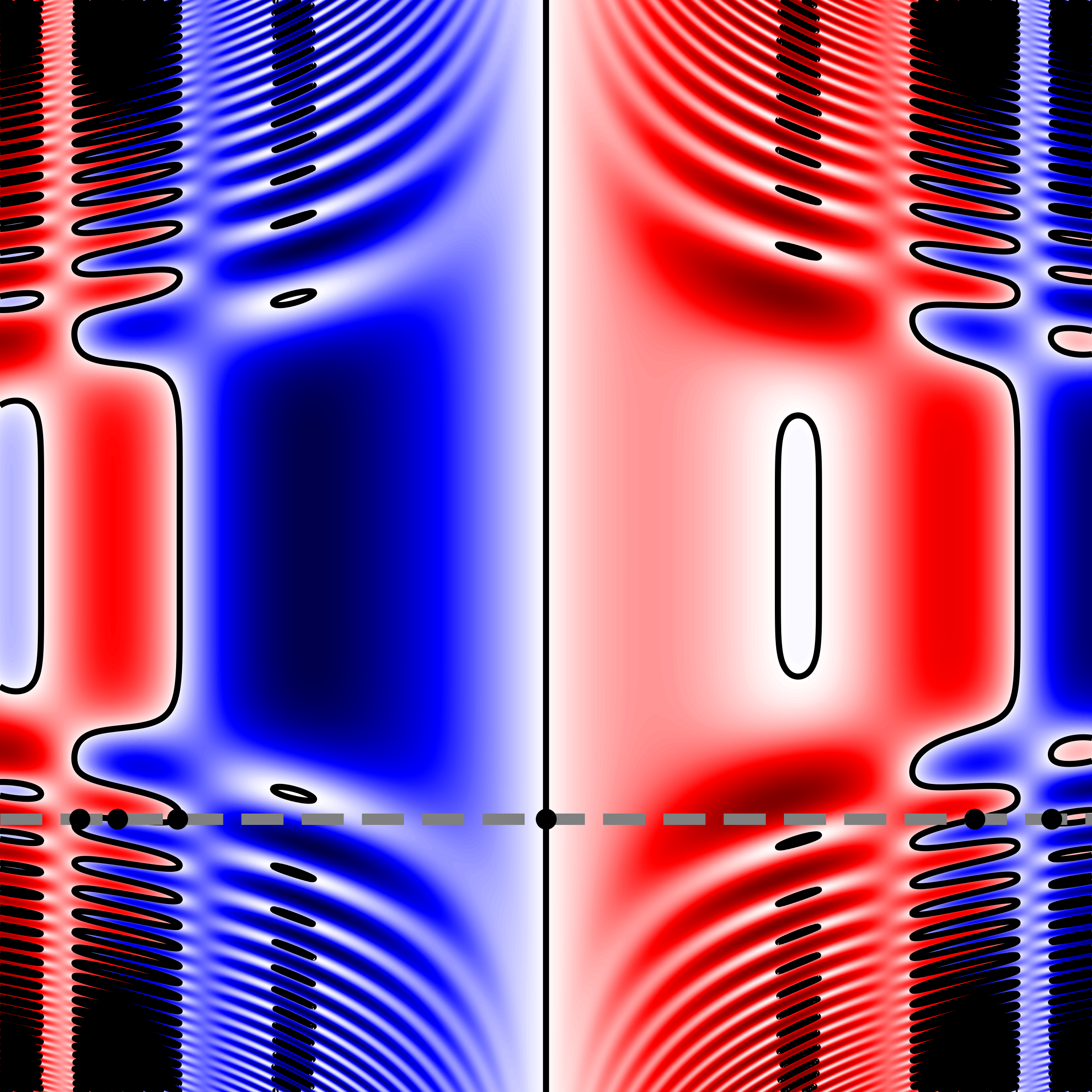};

    \end{axis}

    \node[anchor=north west] at ([xshift=-1.7cm,yshift=0.3cm]axisA.north west) {\labelfontsize (a)};
    \node[anchor=north west] at ([xshift=-1.7cm,yshift=0.3cm]axisB.north west) {\labelfontsize (b)};
    \node[anchor=north west] at ([xshift=-1.7cm,yshift=0.3cm]axisC.north west) {\labelfontsize (c)};

\end{tikzpicture}
  }
  \vspace{-0.5cm}
  \caption{Illustration of the analytic function $\phi$ defined in \cref{eqn:phi} for a two-dimensional input space $\mathcal{F}=\R^2$ and three different parameter-to-recovery maps, along with its zero level set in black line. According to \cref{thm:analytic_varieties_are_small}, this set is compactly slice-null so its intersection (black dots) with a finite dimensional slice, represented by a one-dimensional dashed line, has zero Lebesgue measure within the slice.}
  \label{fig:zer_level_set}
\end{figure}

The function $\phi$ defined in \cref{eqn:phi} is illustrated in \cref{fig:zer_level_set} for two-dimensional spaces $\mathcal{V}=\mathcal{F}=\R^2$ and three different parameter-to-recovery maps. With the notation of \cref{eqn:phi}, we employ canonical basis vectors $e_1^\top$ and $e_2^\top$ as a basis of $\mathcal{V}^*$. Writing $\vect{f}=(f_0,f_1)$, we define
\begin{align*}
  \phi_1(\vect{f}) & = \det\!\begin{bmatrix}
                               \cos\!\big(4 f_0^2\big) & \sin(2f_1+1) \\
                               f_0 f_1                 & 1-f_1^2
                             \end{bmatrix},\quad
  \phi_2(\vect{f})  = \det\!\begin{bmatrix}
                              \sin\!\big(8 f_0^3\big) & 2f_1+1  \\
                              \cos(f_0 f_1)           & 1-f_1^2
                            \end{bmatrix},       \\
  \phi_3(\vect{f}) & = \det\!\begin{bmatrix}
                               \sin(2f_0)      & \sin(2f_0^2) \\
                               \cos(5f_0f_1^4) & 1
                             \end{bmatrix}.
\end{align*}
Here, the dictionary terms of the underlying analytic systems are the columns of the matrices in the determinants. The zero level sets of these functions are displayed in black lines in \cref{fig:zer_level_set}, and have Lebesgue measure zero in $\R$ when intersected with a one-dimensional line.

We now apply \cref{thm:degenerate_sets} to characterize the prevalence of degenerate input sets defined in \cref{def:degenerate_set} for analytic systems of the form of \cref{eqn:analytic_system}.

\begin{corollary}[Prevalence of degenerate sets of analytic systems]\label{cor:analytic_recovery}
  Let $\mcal{F}$ be a Banach space with a countable Schauder basis and suppose \cref{ass:well_posedness} holds. If there exists $f_0 \in \mcal{F} \setminus \Fdeg(\vect{c}^*)$, then $\Fdeg(\vect{c}^*)$ is compactly slice-null in $\mcal{F}$ and has measure zero with respect to any nondegenerate Gaussian measure on Borel sets of $\mcal{F}$.
\end{corollary}
\begin{proof}
  We apply \cref{thm:degenerate_sets} to the mapping $\mcal{F} \to \mcal{B}(\R^N, \mcal{V}) \cong \mcal{V}^N$, $f \mapsto L_{f, \mcal{G}(f)}$, which is analytic by \cref{ass:well_posedness}.
\end{proof}

Finally, the next corollary establishes a zero-one law for the recovery of analytic systems in the form of \cref{eqn:analytic_system}: either coefficient recovery is possible with probability one for input functions sampled from a nondegenerate Gaussian measure, or the degenerate set is the entire space of inputs.

\begin{corollary}[Zero-one law for analytic systems]\label{cor:analytic_recovery_alternative}
  Suppose \cref{ass:well_posedness} holds and $\mcal{F}$ has a countable Schauder basis. Then exactly one of the following is true:
  \begin{enumerate}[label=\roman*), itemsep=0em, leftmargin=*]
    \item If $f$ is sampled from a nondegenerate Gaussian measure on $\mcal{F}$, then $f \notin \Fdeg(\vect{c}^*)$ with probability $1$. \label{recov_1}
    \item $\Fdeg(\vect{c}^*) = \mcal{F}$. \label{recov_2}
  \end{enumerate}
\end{corollary}
\begin{proof}
  If $\mcal{F}\setminus \Fdeg(\vect{c}^*) \neq \emptyset$, then \cref{cor:analytic_recovery} implies \cref{recov_1}. Otherwise, $\Fdeg(\vect{c}^*) = \mcal{F}$, which is \cref{recov_2}.
\end{proof}

\begin{remark}[Recovery up to a constant factor]
  In the situation where $\mcal{D}_0 \equiv 0$ and $\vect{c}^* \neq 0$, one can only hope to determine $\vect{c}^*$ up to a nonzero factor $\lambda\in\R$ from a pair $(f,u)$ as $\sum_{i=1}^N \lambda c^*_i \mcal{D}_i(u, f) = 0$. This is possible when the nullspace $\Null L_{f, u}$ is $1$-dimensional. Then, coefficient recovery up to a constant factor when $\mcal{D}_0 \equiv 0$ is not possible for inputs that belong to the \emph{second-order degenerate set} $\Fdeg^2(\vect{c}^*) = \left\{ f \in \mcal{F} \mid \rank L_{f, \mcal{G}(f)}\leq N-2 \right\}$.
  Analogous zero-one laws hold for this set by the construction given in the proof of \cref{thm:degenerate_sets}.
\end{remark}

\subsection{Multi-shot recovery} \label{sec:multishot}
The same framework covers multi-shot recovery from finitely many input-solution pairs $\{(f_j,u_j)\}_{j=1}^M$ of \cref{eqn:analytic_system} by casting the multi-shot problem as a single-shot problem on product spaces. Given $\vect{u} = (u_1,\ldots,u_M)\in\mcal{U}^M$ and $\vect{f}=(f_1,\ldots,f_M)\in\mcal{F}^M$, we define the stacked parameter-to-output map as the bounded linear operator $L_{\vect{f},\vect{u}}\in\mcal{B}(\R^N,\mcal{V}^M)$ given by
\begin{equation} \label{eq:L_multishot}
  L_{\vect{f},\vect{u}}\vect c = \sum_{i=1}^N c_i\tilde{\mcal{D}}_i(\vect{u},\vect{f})=\sum_{i=1}^N c_i\left(\mcal D_i(u_1,f_1),\ldots,\mcal D_i(u_M,f_M)\right) \in \mcal V^M,
\end{equation}
where $\tilde{\mcal{D}}_i:\mcal{U}^M\times \mcal{F}^M\to\mcal{V}^M$.
One now seeks to identify $\vect{c}^*\in\R^N$ such that
\[ \sum_{i=1}^N c^*_i \tilde{\mcal{D}}_i(\vect{u}, \vect{f}) = \tilde{\mcal{D}}_0(\vect{u},\vect{f}),\]
noting that the unknown coefficient vector is the same as in \cref{eqn:analytic_system}, but the dictionary operators evaluated at multiple input--solution pairs have been ``stacked''. Here, the zero-one mechanism of \cref{cor:analytic_recovery_alternative} applies directly, and we find that either recovery is impossible from $M$ input-output pairs or succeeds with probability $1$ if $\vect{f}$ is sampled from a non-degenerate Gaussian measure on $\mcal{F}^M$. The following proposition shows that adding more input--solution pairs cannot harm recoverability, as expected.

\begin{proposition}[Monotonicity of multi-shot recovery]
  Consider the input-solution pairs $(\vect{f},\vect{u}) = \{(f_j,u_j)\}_{j=1}^M$ and $(\vect{f}',\vect{u}')=\{(f_j,u_j)\}_{j=1}^{M'}$ for $M' \geq M$, and define the associated parameter-to-output recovery maps $L_{\vect{f},\vect{u}}\in\mcal{B}(\R^N,\mcal{V}^{M})$ and $L_{\vect{f}',\vect{u}'}\in\mcal{B}(\R^N,\mcal{V}^{M'})$ as in \cref{eq:L_multishot}. Then,
  \[
    \ker L_{\vect f',\vect u'} \;\subset\; \ker L_{\vect f,\vect u} \quad \text{and}\quad \rank(L_{\vect f',\vect u'}) \geq \rank(L_{\vect f,\vect u}).
  \]
  In particular, if $L_{\vect f,\vect u}$ is injective, then so is $L_{\vect f',\vect u'}$.
\end{proposition}

\begin{proof}
  Let $\vect c\in\ker L_{\vect f',\vect u'}\subset\R^N$, then every coordinate of $L_{\vect{f}',\vect{u}'}\vect c$ vanishes, in particular its first $M$ coordinates:
  \[
    \sum_{i=1}^N c_i\left(\mcal D_i(u_1,f_1),\ldots,\mcal D_i(u_M,f_M)\right)=0,
  \]
  but these are exactly the coordinates of $L_{\vect f,\vect u}\vect c$. Hence $L_{\vect f,\vect u}\vect c = 0$, i.e.\ $\vect c\in\ker L_{\vect f,\vect u}$. The rank inequality follows from the rank--nullity theorem applied to both operators on the same domain $\R^N$.
\end{proof}

\subsection{The one-shot identification algorithm}

\Cref{sec:genericity} characterizes the injectivity of the map $L_{f,u}:\R^N\to\mcal{V}$. We introduce a computable weak-form method to test the finite-dimensional map $\Psi L_{f,u}$, where $\Psi:\mcal{V}\to\R^{r_0}$ is defined by bounded linear functionals $\psi_1,\ldots,\psi_{r_0}\in\mcal{V}^*$, with $r_0\geq N$. The two notions agree when $\Psi$ is injective on $\Range(L_{f,u})$.

\begin{algorithm}[htbp]
  \caption{Weak Identification of Analytic Systems (WIAS)}
  \label{alg:weak_identification}
  \begin{algorithmic}[1]
    \Require Source term $f$, observed solution $u$, dictionary $\{\mathcal{D}_i\}_{i=0}^N$, dual space family $\{\psi_i\}_{i=1}^{r_0}\subset \mcal{V}^*$
    \State Compute the $r_0\times N$ matrix $G_{i,j} = \psi_i(\mathcal{D}_j(u,f))$
    \State Compute the right-hand side $b_i = \psi_i(\mathcal{D}_0(u,f))$, $1\leq i\leq r_0$
    \State Compute a column-pivoted QR factorization $GP = QR$
    \State Determine the numerical rank $r$ from the diagonal of $R$
    \If{$r<N$}
    \State \Return non-identifiability from this input and test set
    \EndIf
    \State Solve $R_{1:N,1:N} \vect{z} = Q_{:,1:N}^\top \vect{b}$
    \State Recover the coefficient vector by undoing the permutation: $\vect{c}=P\vect{z}$
    \State \Return $\vect{c}$
  \end{algorithmic}
\end{algorithm}

Similarly to Weak SINDy~\cite{messenger2021weak,messenger2021weakgalerking}, we select test functionals $\psi_1,\ldots,\psi_{r_0}\in \mcal{V}^*$ with $r_0\geq N$. When $\mcal{V}$ is a Hilbert space, such functionals may be represented by inner products via the Riesz representation theorem. If the space has been discretized with finite elements, one may instead choose a basis for the finite-dimensional dual space. We then assemble an overdetermined linear system $G \vect{c} = b$, where $G_{i,j} = \psi_i(\mathcal{D}_j(u,f))$ and $b_i = \psi_i(\mathcal{D}_0(u,f))$ for $1\leq i\leq r_{0}$ and $1\leq j\leq N$. In practice, we employ a rank-revealing column-pivoted QR factorization to solve this system. The identification of the numerical rank of $G$ allows us to determine a posteriori whether the source term $f$ lies in the degenerate set $\Fdeg(\vect{c}^*)$ and guarantees the correctness of the method according to \cref{cor:tested_recovery}. The procedure, which we call Weak Identification of Analytic Systems (WIAS), is described in \cref{alg:weak_identification}.

\begin{corollary}[Zero-one law for tested recovery]\label{cor:tested_recovery}
  Let $\Psi:\mcal{V}\to\R^{r_0}$ be a bounded linear map. Under the hypotheses of \cref{thm:degenerate_sets}, define
  \[
    \mcal{F}_{{\rm deg},\Psi}=\{f\in\mcal{F}\mid \rank(\Psi L_f)\leq N-1\}.
  \]
  If there exists $f_0\in\mcal{F}$ such that $\rank(\Psi L_{f_0})=N$, then $\mcal{F}_{{\rm deg},\Psi}$ is compactly slice-null and has measure zero with respect to any nondegenerate Gaussian measure on Borel sets of $\mcal{F}$.
\end{corollary}
\begin{proof}
  Apply \cref{thm:degenerate_sets} to the analytic mapping $f\mapsto \Psi L_f$ with codomain $\mcal{B}(\R^N,\R^{r_0})$.
\end{proof}

Similarly to \cref{cor:analytic_recovery}, after applying \cref{cor:tested_recovery} to the mapping $f\mapsto L_{f,\mathcal{G}(f)}$, we find that $\Fdeg(\vect{c}^*)\subset  \mcal{F}_{{\rm deg},\Psi}(\vect{c}^*)$ as $\rank(L_{f,\mathcal{G}(f)})<N$ implies that $\Psi L_{f,\mathcal{G}(f)}$ is not full rank. Hence, the existence of a single input $f_0\in \Fdeg(\vect{c}^*)$ implies that the degenerate set $\Fdeg(\vect{c}^*)$ is compactly slice-null and has measure zero with respect to any nondegenerate Gaussian measure on Borel sets of $\mcal{F}$. This guarantees correctness of \cref{alg:weak_identification} but does not ensure that failing inputs lie in the degenerate set $\Fdeg(\vect{c}^*)$ as the test functionals may be poorly selected. In practice, one should choose a sufficiently large number of test functionals to ensure that $\Psi$ is injective on $\Range(L_{f,\mathcal{G}(f)})$ for all $f\in\mcal{F}$. If $\mathcal{V}$ is a Hilbert space, one may select $\Psi = L_{f,\mathcal{G}(f)}^*$ with the functionals $\psi_1=\langle\cdot, \mathcal{D}_1(u,f)\rangle,\ldots,\psi_{N}=\langle\cdot, \mathcal{D}_{N}(u,f)\rangle$ to obtain $\Fdeg(\vect{c}^*) = \mcal{F}_{{\rm deg},\Psi}(\vect{c}^*)$. In this case, failing inputs are guaranteed to lie in the degenerate set $\Fdeg(\vect{c}^*)$.

\section{Matrix Recovery}
\label{sec:matrix_recovery}

We show how the framework from \Cref{sec:analytic_systems} can be applied to recover matrices in linearly parameterized families using matrix-vector products \cite{halikias2022matrix}.
A linearly parameterized family of matrices is an $N$-dimensional subspace $\mcal{A} \subset \R^{m\times n}$ of $m\times n$ matrices with a basis $\{ A_i \}_{i=1}^{N}$.
Matrices in this family are parametrized by the linear map $A: \R^N \to \mcal{A}$ defined by
\[
  A(\vect{c}) = \sum_{i=1}^{N} c_i A_i, \quad \vect{c}\in\R^N.
\]
Examples of linearly parametrized families include diagonal, circulant, Toeplitz, Hankel, and matrices with fixed sparsity pattern. Related works in numerical linear algebra study the near-optimal approximation problem for the class of matrices with a fixed sparsity pattern~\cite{amsel2026fixed, musco2026nearly}. More recently, a different approach based on finite discretizations of a linearly parametrized family was proposed in~\cite{amsel2026query}.

In this section, we derive results about degenerate sets in two types of matrix recovery problems for an unknown element $A(\vect{c}^*)$ from the linear family accessible through matrix-vector products. For the first type of matrix recovery problem, we select $q_R$ ``right'' and $q_L$ ``left'' input vectors $f_1, \ldots, f_{q_R} \in \R^n$ and $g_1, \ldots, g_{q_L} \in \R^{m}$, and acquire the corresponding outputs $u_{1}, \ldots, u_{q_R} \in \R^m$ and $v_{1}, \ldots, v_{q_L} \in \R^n$ by applying the unknown matrix $A(\vect{c}^*)$ and its transpose to the inputs, that is: recover $\vect{c}^* \in \R^N$ from the input-output pairs $(f_j, u_j)$ and $(g_k, v_k)$ such that
\[
  A(\vect{c}^*) f_j = u_j
  \quad \text{and} \quad
  A(\vect{c}^*)^\top g_k = v_k,\qquad 1\leq j \leq q_R, \, 1\leq k \leq q_L.
\]
We ask whether these input-output pairs are sufficient to recover the unknown coefficients $\vect{c}^*$. To cast this problem within the framework of \Cref{sec:multishot}, we observe that the combined inputs $\vect{f} = (f_1, \ldots, f_{q_R}, g_1, \ldots, g_{q_L}) \in \mcal{F} = \R^{q_R n + q_L m}$ and outputs $\vect{u} = (u_1, \ldots, u_{q_R}, v_1, \ldots, v_{q_L}) \in \mcal{U} = \R^{q_R m + q_L n}$ satisfy
\begin{equation}
  \mcal{L}(\vect{u}, \vect{f})
  \coloneqq
  \sum_{i=1}^{N} c_i^*
  \underbrace{
    \begin{bmatrix}
      A_i &        &     &          &        &            \\
          & \ddots &     &          & 0      &            \\
          &        & A_i &          &        &            \\
          &        &     & A_i^\top &        &            \\
          & 0      &     &          & \ddots &            \\
          &        &     &          &        & A_i^{\top}
    \end{bmatrix}
    \begin{bmatrix}
      f_1     \\
      \vdots  \\
      f_{q_R} \\
      g_1     \\
      \vdots  \\
      g_{q_L}
    \end{bmatrix}
  }_{\mcal{D}_i(\vect{u},\vect{f})}
  = \underbrace{
    \begin{bmatrix}
      u_1     \\
      \vdots  \\
      u_{q_R} \\
      v_1     \\
      \vdots  \\
      v_{q_L}
    \end{bmatrix}
  }_{\mcal{D}_0(\vect{u},\vect{f})}.
\end{equation}
In this case, we seek to establish whether the linear recovery operator
\[
  L_{\vect{f}}:\vect{c}\mapsto
  \begin{bmatrix}
    \sum_{i=1}^{N}c_i A_i f_1      \\
    \vdots                         \\
    \sum_{i=1}^{N}c_i A_i f_{q_R}  \\
    \sum_{i=1}^{N}c_i A_i^\top g_1 \\
    \vdots                         \\
    \sum_{i=1}^{N}c_i A_i^\top g_{q_L}
  \end{bmatrix}
\]
is injective for generic choices of $\vect{f} = (f_1, \ldots, f_{q_R}, g_1, \ldots, g_{q_L})\in \R^{q_R n + q_L m}$. Following \cref{def:degenerate_set}, we introduce the degenerate set $\Fdeg(\vect{c}^*)$ of inputs $\vect{f}$ for which $L_{\vect{f}}$ is non-injective. The following result recovers the main result of \cite{otto2023note} and establishes a zero-one law for the recovery of the unknown coefficients $\vect{c}^*$, hence the matrix $A$ within the linearly parameterized family $\vspan\{A_1,\ldots,A_N\}$, from the above input-output pairs: either $\vect{c}^*$ can be recovered with probability one for inputs sampled from a nondegenerate Gaussian measure or it cannot be recovered for any input.

\begin{corollary}[Zero-one law for matrix recovery]\label{cor:matrix_recovery}
  Exactly one of the following is true:
  \begin{enumerate}[label=\roman*), itemsep=0em, leftmargin=*]
    \item If the vector $\vect{f} = (f_1, \ldots, f_{q_R}, g_1, \ldots, g_{q_L}) \in \R^{q_R n + q_L m}$ is sampled from a nondegenerate Gaussian measure on $\mcal{F}$, then $\vect{f} \notin \Fdeg(\vect{c}^*)$ with probability $1$. \label{recov_1_matrix}
    \item $\Fdeg(\vect{c}^*) = \mcal{F}$. \label{recov_2_matrix}
  \end{enumerate}
\end{corollary}

\begin{proof}
  This is clearly in the form of \cref{eqn:analytic_system} with linear, hence analytic dictionary functions and solution operator $\mcal{G} : \vect{f} \mapsto \vect{u}$ defined by $\mcal{G}(\vect{f})=\big((I_{q_R}\otimes A(\vect{c}^*))\oplus (I_{q_L}\otimes A(\vect{c}^*)^\top)\big)\vect{f}$.
  The degenerate set $\Fdeg(\vect{c}^*)$ consists of the inputs $\vect{f}$ for which there are other coefficients $\vect{c} \neq \vect{c}^*$ satisfying the above equation with the same observed output $\vect{u} = \mcal{G}(\vect{f})$.
\end{proof}

In the second type of matrix recovery problem, we assume that the unknown matrix is accessible through matrix-vector products with its inverse.
Specifically, we restrict to the square case $m=n$, assume that $A(\vect{c}^*)$ is invertible, and after selecting $q_R$ and $q_L$ right and left input vectors $f_1, \ldots, f_{q_R} \in \R^n$ and $g_1, \ldots, g_{q_L} \in \R^{n}$,
we observe outputs $u_{1}, \ldots, u_{q_R} \in \R^n$ and $v_{1}, \ldots, v_{q_L} \in \R^n$ defined by
\[
  u_j = A(\vect{c}^*)^{-1} f_j
  \quad \text{and} \quad
  v_k = A(\vect{c}^*)^{-\top} g_k.
\]
The same combined inputs $\vect{f}$ and combined outputs $\vect{u}$ satisfy
\begin{equation}
  \mcal{L}(\vect{u}, \vect{f})
  \coloneqq
  \sum_{i=1}^{N} c_{i}^*
  \underbrace{
    \begin{bmatrix}
      A_i &        &     &          &        &          \\
          & \ddots &     &          &        &          \\
          &        & A_i &          &        &          \\
          &        &     & A_i^\top &        &          \\
          &        &     &          & \ddots &          \\
          &        &     &          &        & A_i^\top
    \end{bmatrix}
    \begin{bmatrix}
      u_1     \\
      \vdots  \\
      u_{q_R} \\
      v_1     \\
      \vdots  \\
      v_{q_L}
    \end{bmatrix}
  }_{\mcal{D}_i(\vect{u},\vect{f})}
  =
  \underbrace{
    \begin{bmatrix}
      f_1     \\
      \vdots  \\
      f_{q_R} \\
      g_1     \\
      \vdots  \\
      g_{q_L}
    \end{bmatrix}
  }_{\mcal{D}_0(\vect{u},\vect{f})},
\end{equation}
and the linear solution operator $\mcal{G}: \vect{f} \mapsto \vect{u}$ is analytic.
In this case, \cref{cor:analytic_recovery} produces the same corollary as \cref{cor:matrix_recovery} for the respective degenerate sets.

\begin{remark}[Recovery from sketching $A$ and its inverse]
  Accessing matrix-vector products with all four operators $A(\vect{c}^*)$, $A(\vect{c}^*)^\top$, $A(\vect{c}^*)^{-1}$, and $A(\vect{c}^*)^{-\top}$ can be treated analogously, and \cref{cor:analytic_recovery} produces the same corollary for the corresponding degenerate set.
\end{remark}

\begin{remark}[No adjoint access]
  If adjoint products are unavailable, one may take $q_L=0$ and apply the same zero-one law to right matrix-vector products only. Recovery then depends on whether the corresponding right-product recovery map has full column rank for at least one batch of inputs.
\end{remark}

\section{Analytic continuous-time dynamical systems} \label{sec:dyn_systems}

We apply the framework from \Cref{sec:analytic_systems} to show that if the coefficients describing an analytic continuous-time dynamical system are recoverable using the trajectory produced by a certain input, then they are recoverable using trajectories generated by all inputs except for a set that is compactly slice-null and has zero nondegenerate Gaussian measure.
Specifically, we consider a dynamical system described by a first-order ordinary differential equation (ODE) in the sense of Carath\'{e}odory:
\begin{equation}\label{eqn:analytic_ODE}
  \frac{\rm{d}}{\rm{d} t} u(t) = F(t, u(t), f(t)) \coloneqq \sum_{i = 1}^{N} c^*_i D_i(t, u(t), f(t)),
  \qquad u(a) = f_a,
\end{equation}
for an absolutely continuous function $u \in \mcal{U} = \mcal{V} \coloneqq \AC([a,b];\R^n)$ on a bounded interval $[a,b] \subset \R$, $a < b$, with inputs $(f_a, f)$ in a closed subspace $\mcal{F} \subset \R^n \times L^{\infty}([a,b]; \R^m)$. Here, we observe the system dynamics on a bounded time interval $[a,b]$.
The zero-one results below apply to those input subspaces that additionally possess a countable Schauder basis, such as $\R^n\times C([a,b];\R^m)$.
We assume the dictionary functions $D_i : [a, b] \times \R^n \times \R^m \to \R^n$ are known, allowing us to express \cref{eqn:analytic_ODE} in the form of \cref{eqn:analytic_system} with dictionary operators $\mcal{D}_i : \mcal{U} \times \mcal{F} \to \mcal{U}$ defined by
\begin{equation}\label{eqn:analytic_ODE_dict_ops}
  \mcal{D}_0[u,(f_a, f)](t) = u(t) - f_a,
  \quad
  \mcal{D}_i[u,(f_a, f)](t) = \int_{a}^{t} D_i(s, u(s), f(s)) \d s,
\end{equation}
for $t\in[a,b]$ and $1\leq i\leq N$.
With the notations of \cref{eqn:analytic_system}, we seek to recover the operator $\mathcal{L}:\mathcal{U}\times \mathcal{F}\to\mathcal{U}$ defined for $u\in\mathcal{U}, (f_a,f)\in\mathcal{F}$ and $t\in[a,b]$ by
\begin{equation} \label{eq:operator_L_dyn}
  \mathcal{L}[u,(f_a,f)](t) = \int_{a}^{t} F(s, u(s), f(s)) \d s = \sum_{i=1}^N c^*_i\int_{a}^{t} D_i(s, u(s), f(s)) \d s.
\end{equation}
We make the following assumptions and show that they ensure the framework from \Cref{sec:analytic_systems} can be applied.

\begin{assumption}[Regularity conditions for dynamical systems]\label{ass:integrable_analytic_ODE}
  \leavevmode
  \begin{enumerate}[label=\roman*), itemsep=0em, leftmargin=*]
    \item For every $(f_a, f) \in \mcal{F}$, \cref{eqn:analytic_ODE} has a solution $u \in \mcal{U}$. \label{asmpn:ODE_existence}
    \item For every $(x,y)\in\R^n\times \R^m$, $t\mapsto D_i(t,x,y)$ is measurable and $(x,y) \mapsto D_i(t, x, y)$ is analytic for almost every $t \in [a,b]$. \label{asmpn:ODE_dictionary_regularity}
    \item \label{item:assum_ODE_3} For every compact subset $\mcal{C} \subset \R^n \times \R^m$, there exists a radius $r > 0$ and a constant $M_\mcal{C}$ such that for almost every $t \in [a,b]$, each function $(x,y) \mapsto D_i(t,x,y)$ admits a holomorphic extension to the closed complex neighborhood
          \[
            B_{\C}(\mcal{C}, r) = \left\{ (w,z) \in \C^{n} \times \C^{m} \mid \exists (x,y) \in \mcal{C} \ \mbox{s.t.} \ \| (w,z) - (x,y) \| \leq r \right\},
          \]
          that is bounded by $\sup_{(w,z) \in B_{\C}(\mcal{C}, r)} \left\| D_i(t, w, z) \right\| \leq M_\mcal{C}$.
  \end{enumerate}
\end{assumption}

Note that if the dictionary functions satisfy \cref{ass:integrable_analytic_ODE} \cref{item:assum_ODE_3}, then the function $F$ satisfies the same property by linearity.

\begin{example}[Time-invariant systems]
  \cref{ass:integrable_analytic_ODE} \cref{asmpn:ODE_dictionary_regularity} and \cref{item:assum_ODE_3} are satisfied if the dynamical system is time-invariant, i.e, of the form
  \[\frac{\rm{d}}{\rm{d} t} u(t) = F(u(t), f(t)),\]
  and the dictionary functions $D_i$ only depend on $u(t)$ and $f(t)$ and are analytic with respect to these arguments. This includes standard polynomial dictionaries used for systems such as the Lorenz equations and Duffing oscillator.
\end{example}

The following lemmas show that under \cref{ass:integrable_analytic_ODE}, the dictionary operators defined by \cref{eqn:analytic_ODE_dict_ops} are analytic, the solution of \cref{eqn:analytic_ODE} is unique, and it is given by an analytic operator $\mcal{G} : (f_a, f) \mapsto u$.
Thus, \cref{ass:well_posedness} holds as a consequence of \cref{ass:integrable_analytic_ODE}, and
\cref{cor:analytic_recovery,cor:analytic_recovery_alternative} apply directly to coefficient recovery for \cref{eqn:analytic_ODE}.

\begin{lemma}[Well-posedness and continuity of $\mcal{G}$] \label{lem:continuity_of_ODE_solution}
  Under \cref{ass:integrable_analytic_ODE}, the solution $u\in \mcal{U}$ to \cref{eqn:analytic_ODE} is unique for each $(f_a, f) \in \mcal{F}$, and the resulting solution operator $\mcal{G}: \mcal{F} \to \mcal{U}$ defined by $\mathcal{G} (f_a, f) = u$
  is continuous.
\end{lemma}
\begin{proof}
  The proof follows an approach similar to \cite[Thm.~3.5]{khalil2002nonlinear}.
  Let $(f_{a,*}, f_*)\in\mcal{F}$ and $u_* \in \mcal{U}$ be a corresponding solution to \cref{eqn:analytic_ODE} guaranteed by \cref{ass:integrable_analytic_ODE} \cref{asmpn:ODE_existence}.
  Choose $0 < \varepsilon < 1$ and define the closed tube around the graph of the solution by
  \[
    \mcal{U}_{\varepsilon} = \{ (t, x) \in [a,b]\times \R^n \ : \ \| x - u_*(t) \| \leq \varepsilon \}.
  \]
  Consider the compact subset $\mcal{C} = \mcal{C}_u \times \mcal{C}_f \subset \R^{n}\times\R^m$ where
  \[
    \mcal{C}_u = \{ x \in \R^n \ : \ \| x \| \leq \sup_{t\in [a,b]} \|u_*(t)\|+1 \}, \,\,
    \mcal{C}_f = \{ y \in \R^m \ : \ \| y \| \leq \esssup_{t\in [a,b]} \|f_*(t)\|+1 \}.
  \]
  Note that $\sup_{t\in [a,b]} \|u_*(t)\|$ exists by the extreme value theorem since $u$ is absolutely continuous. With this compact subset and the holomorphic extension of $(x,y) \mapsto F(t,x,y)$, denoted $\check{F}$, to a neighborhood $B_{\C}(\mcal{C}, r)$ with bound $M=\sum_{i=1}^N|c^*_i|M_{\mcal{C}}$ provided by \cref{ass:integrable_analytic_ODE} \cref{item:assum_ODE_3}, applying Cauchy's integral formula to the function $\zeta\mapsto \check{F}(t, x + v\zeta, y + w \zeta)$ with contour $\gamma_r : \theta \mapsto r e^{i \theta}$, for $\theta\in[0,2\pi)$ gives
  \[
    \left.\frac{\d}{\d \zeta} \check{F}(t, x + v\zeta, y + w \zeta)\right\vert_{\zeta = 0} = \frac{1}{2\pi i} \oint_{\gamma_r} \frac{\check{F}(t, x + v\zeta, y + w \zeta)}{\zeta^2} \d \zeta,
  \]
  for almost every $t \in [a,b]$, every $(x,y) \in \mcal{C}$ and $(v,w) \in \C^{n}\times \C^{m}$ with $\| (v,w) \| \leq 1$.
  In particular, this means that
  \begin{equation}\label{eqn:uniform_bound_on_derivatives_of_F}
    \begin{aligned}
      \left\Vert \frac{\partial F}{\partial x}(t, x, y) \right\Vert
                                                                    & = \sup_{\substack{v\in \R^n \\ \| v \|\leq 1}} \left\| \left.\frac{\d}{\d \zeta} \check{F}(t, x + v\zeta, y)\right\vert_{\zeta = 0} \right\|
      \leq \frac{M}{r} =: L_u,                                                                                                                                                                                     \\
      \left\Vert \frac{\partial F}{\partial y}(t, x, y) \right\Vert & = \sup_{\substack{w\in \R^m \\ \| w \|\leq 1}} \left\| \left.\frac{\d}{\d \zeta} \check{F}(t, x, y + w\zeta)\right\vert_{\zeta = 0} \right\|
      \leq \frac{M}{r} =: L_f.
    \end{aligned}
  \end{equation}
  Let $\varepsilon>0$, $\delta_{\varepsilon} = \min\{1, \ \frac{\varepsilon}{1 + (b-a)L_f} e^{-L_u(b-a)} \}$,
  and choose any $(f_a, f) \in \mcal{F}$ such that
  \[
    \| (f_a, f) - (f_{a,*}, f_*) \|
    = \| f_a - f_{a,*} \|_{\R^n} + \| f - f_* \|_{L^{\infty}([a,b];\R^m)}
    < \delta_{\varepsilon}.
  \]
  Since this implies that $\| f(t)\| - \|f_*(t) \|\leq \| f(t) - f_*(t) \| < 1$ for almost every $t \in [a,b]$, it follows that
  for any such $t$ we have $f_*(t), f(t) \in \mcal{C}_f$ and for every $x \in \mcal{C}_u$, \cref{eqn:uniform_bound_on_derivatives_of_F} implies
  \begin{equation}\label{eqn:lem_f_Lipschitz_estimate}
    \| F(t, x, f(t)) - F(t, x, f_*(t)) \| \leq L_f \| f(t) - f_*(t) \| \leq L_f \| (f_a, f) - (f_{a,*}, f_*) \|.
  \end{equation}
  Let $u$ be a solution to \cref{eqn:analytic_ODE} with input $(f_a, f)$ and $\tau \coloneqq \sup \{ s \in [a,b] \ : \ (t, u(t)) \in \mcal{U}_{\varepsilon} \ \text{for every } t\in[a,s] \}$. This set is nonempty because $\|f_a-f_{a,*}\|<\varepsilon$, and the supremum is attained by continuity of $u$ and closedness of the tube.
  For every $t \in [a, \tau]$ we have $(t, u(t)) \in \mcal{U}_{\varepsilon}$, which implies that $u(t) \in \mcal{C}_u$. Hence, for almost every $t \in [a, \tau]$ we have
  \begin{equation}\label{eqn:lem_u_Lipschitz_estimate}
    \| F(t, u(t), f(t)) - F(t, u_*(t), f(t)) \| \leq L_u \| u(t) - u_*(t) \|,
  \end{equation}
  by \cref{eqn:uniform_bound_on_derivatives_of_F} since $f\in\mathcal{C}_f$. Let $t\in[a,\tau]$, integrating \cref{eqn:analytic_ODE} between $a$ and $t$ and using the Lipschitz estimates from \cref{eqn:lem_f_Lipschitz_estimate,eqn:lem_u_Lipschitz_estimate} yield
  \begin{align*}
    \| u(t) - u_*(t) \|
     & \leq \| f_a - f_{a,*} \| + \int_{a}^{t} \| F(s, u(s), f(s)) - F(s, u_*(s), f_*(s)) \| \d s                       \\
     & \leq \| f_a - f_{a,*} \|
    \begin{aligned}[t]
       & +\int_{a}^{t} \| F(s, u(s), f(s)) - F(s, u_*(s), f(s)) \| \,\d s      \\
       & + \int_{a}^{t} \| F(s, u_*(s), f(s)) - F(s, u_*(s), f_*(s)) \| \,\d s
    \end{aligned}                                          \\
     & \leq \| f_a - f_{a,*} \| + \int_{a}^{t} L_u \| u(s) - u_*(s) \| \d s + \int_{a}^{t} L_f \| f(s) - f_*(s) \| \d s \\
     & \leq (1 + (b-a)L_f) \| (f_a, f) - (f_{a,*}, f_*) \| + \int_{a}^{t} L_u \| u(s) - u_*(s) \| \d s.
  \end{align*}
  By the Gr\"{o}nwall-Bellman inequality~\cite[Lem.~A.1]{khalil2002nonlinear} and our choice of $\delta_{\varepsilon}$, we obtain
  \[
    \| u(t) - u_*(t) \| \leq (1 + (b-a)L_f) \| (f_a, f) - (f_{a,*}, f_*) \| e^{L_u (t-a)} < \varepsilon, \quad \forall t \in [a, \tau].
  \]
  If $\tau < b$ then by continuity of $u$ and $u_*$, the above inequality continues to hold on $[a,\tau']$ for some $\tau'>\tau$, meaning $(t,u(t)) \in \mcal{U}_{\varepsilon}$ for all $t \in [a,\tau']$. Since this contradicts the definition of $\tau$, it follows that $\tau = b$ and $\| u(t) - u_*(t) \| < \varepsilon$ for all $t \in [a,b]$.
  This proves that \cref{eqn:analytic_ODE} has unique solutions on $[a,b]$ and that the resulting mapping $\mcal{G}$ is continuous from $\mcal{F}$ into $C([a,b]; \R^n)$.

  To prove continuity of $\mcal{G}$ into $\mcal{U} = \AC([a,b]; \R^n)$, which is a Banach space equivalent to the Sobolev space $W^{1,1}([a,b];\R^n)$ \cite[Rem.~8, p.~206]{brezis2010functional}, we observe that
  \begin{align*}
    \int_{a}^{b} \left\| \frac{\d}{\d t} u(t) - \frac{\d}{\d t} u_*(t) \right\| \d t
     & = \int_{a}^{b} \left\| F(t, u(t), f(t)) - F(t, u_*(t), f_*(t)) \right\| \d t                                                                                                                  \\
     & \hspace{-3.8cm} \leq \int_{a}^{b} \!\!\! \left(\left\| F(t, u(t), f(t)) \!-\! F(t, u_*(t), f(t)) \right\|\! +\! \left\| F(t, u_*(t), f(t)) \!-\! F(t, u_*(t), f_*(t)) \right\| \right)\! \d t \\
     & \hspace{-3.8cm} \leq \int_{a}^{b} \left(L_u \left\| u(t) - u_*(t) \right\| + L_f \left\| f(t) - f_*(t) \right\| \right) \d t                                                                  \\
     & \hspace{-3.8cm} < (b-a) \left(L_u \| \mcal{G}(f_a, f) - \mcal{G}(f_{a,*}, f_*) \|_{C([a,b];\R^n)} + L_f \| f - f_* \|_{L^{\infty}([a,b];\R^m)}\right),
    % \to 0 \quad \mbox{as} \quad \| (f_a, f) - (f_{a,*}, f_*) \| \to 0
  \end{align*}
  where the second inequality comes from \cref{eqn:uniform_bound_on_derivatives_of_F}. Following the continuity of $\mcal{G}$ from $\mcal{F}$ to $C([a,b];\R^n)$, we obtain
  \[
    \| \mcal{G}(f_a, f) - \mcal{G}(f_{a,*}, f_*) \|_{W^{1,1}([a,b];\R^n)} \longrightarrow 0 \quad \mbox{as} \quad \| (f_a, f) - (f_{a,*}, f_*) \|_{\mcal{F}} \longrightarrow 0,
  \]
  which proves that $\mcal{G} : \mcal{F} \to \mcal{U}$ is continuous.
\end{proof}

The following is a technical lemma establishing analyticity of the dictionary operators and providing explicit expressions for their first-order Fr\'{e}chet derivatives.

\begin{lemma}[Analyticity of dictionary operators]\label{lem:analytic_Nemytskii_operator}
  Under \cref{ass:integrable_analytic_ODE}, each operator $\mcal{D}_i : \mcal{U} \times \mcal{F} \to \mcal{U}$ defined in \eqref{eqn:analytic_ODE_dict_ops} is analytic and the action of its Fr\'{e}chet derivative $\D \mcal{D}_i[u, (f_a, f)] : \mcal{U} \times \mcal{F} \to \mcal{U}$ at each $(u, (f_a, f))\in \mcal{U}\times \mcal{F}$ on $(v, (w_a, w)) \in \mcal{U}\times \mcal{F}$ is given for $t\in[a,b]$ and $1\leq i\leq N$ by
  \[
    \D \mcal{D}_i[u, (f_a, f)](v, (w_a, w))(t)
    \!= \!\!\int_{a}^{t}\! \left[
      \frac{\partial D_i}{\partial x}(s, u(s), f(s)) v(s)
      \!+\! \frac{\partial D_i}{\partial y}(s, u(s), f(s)) w(s)\!
      \right]\!\! \!\d s,
  \]
  where $\partial_x D_i(s, u(s), f(s))$ and $\partial_y D_i(s, u(s), f(s))$ are uniformly bounded for almost all $s \in [a,b]$, and  $\D \mcal{D}_0[u, (f_a, f)](v, (w_a, w)) = v - w_a$.
\end{lemma}
\begin{proof}
  The result for $\mcal{D}_0$ is obvious because it is a bounded linear operator, and thus analytic.
  Throughout the rest of the proof, we fix index $1\leq i\leq N$ and consider the operator $\mcal{D} = \mcal{D}_i$ with dictionary function $D = D_i$.
  Let $(u_*, (f_{a,*}, f_*)) \in \mcal{U}\times \mcal{F}$.
  To show that $\mcal{D}$ is analytic near this point, we consider a holomorphic extension of $\mcal{D}$ to a neighborhood of $(u_*, (f_{a,*}, f_*))$ in the complex Banach space $\check{\mcal{U}} \times \check{\mcal{F}}$ with $\check{\mcal{U}} = \AC([a,b];\C^n)$ and $\check{\mcal{F}} = \C^n \times L^{\infty}([a,b];\C^m)$ defined as follows.
  With the compact set $\mcal{C} = \mcal{C}_u \times \mcal{C}_f \subset \R^{n}\times\R^m$, where
  \[
    \mcal{C}_u = \left\{ x \in \R^n \mid \| x \| \leq \!\sup_{t\in [a,b]} \|u_*(t)\|\!+1\!\! \right\},\,
    \mcal{C}_f = \left\{ y \in \R^m \mid \| y \| \leq \esssup_{t\in [a,b]} \|f_*(t)\|\!+1\! \!\right\},
  \]
  \cref{ass:integrable_analytic_ODE} \cref{item:assum_ODE_3} ensures that $D$ admits a holomorphic extension $\check{D}: [a,b] \times B_{\C}(\mcal{C},r) \to \C^n$ with magnitude bounded by $M_{\mcal{C}}$ and $0 < r < 1$.
  We introduce the following open neighborhoods of $u_*$ and $(f_{a,*}, f_*)$:
  \begin{align*}
    \mcal{B}_{\C}(u_*, r/4)            & =
    \left\{ u \in \check{\mcal{U}} \ : \ \| u - u_* \| < r/4 \right\}, \\
    \mcal{B}_{\C}((f_{a,*}, f_*), r/4) & =
    \left\{ (f_a, f) \in \check{\mcal{F}} \ : \ \| f_a - f_{a,*} \| + \| f - f_* \| < r/4 \right\},
  \end{align*}
  and define the holomorphic extension $\check{\mcal{D}}: \mcal{B}_{\C}(u_*, r/4) \times \mcal{B}_{\C}((f_{a,*}, f_*), r/4) \subset \check{\mcal{U}} \times \check{\mcal{F}} \to \check{\mcal{U}}$ by
  \[
    \check{\mcal{D}}[u, (f_a, f)](t) = \int_{a}^{t} \check{D}(s, u(s), f(s)) \d s, \qquad t \in [a,b].
  \]
  This is well-defined because $u \in \mcal{B}_{\C}(u_*, r/4)$ and $(f_a, f) \in \mcal{B}_{\C}((f_{a,*}, f_*), r/4)$ implies that for almost every $t \in [a,b]$ we have,
  \[
    \| f_a - f_{a,*} \| + \| f(t) - f_*(t) \| + \| u(t) - u_*(t) \| < r/2,
  \]
  meaning $(u(t), f(t)) \in B_{\C}(\mcal{C}, r)$, where $\check{D}(t, \cdot, \cdot)$ is defined.
  We prove that $\check{\mcal{D}}$ is holomorphic using \cite[Thm.~4]{bochnak1970analytic} by establishing that $\check{\mcal{D}}$ is locally bounded and holomorphic on $1$-dimensional affine subspaces.
  Local boundedness is obvious because
  \[
    \left\| \check{\mcal{D}}(u,(f_a, f)) \right\|_{W^{1,1}([a,b];\C^n)}
    \leq \frac{1}{2}(b-a)^2 M_{\mcal{C}}
    + (b-a)M_{\mcal{C}}.
  \]
  Let $(u, (f_a, f)) \in \mcal{B}_{\C}(u_*, r/4) \times \mcal{B}_{\C}((f_{a,*}, f_*), r/4)$ and $(v, (w_a, w)) \in \check{\mcal{U}} \times \check{\mcal{F}}$, and consider an affine subspace $\Lambda = (u, (f_a, f)) + \vspan\{ (v, (w_a, w)) \}$.
  To show that $\check{\mcal{D}}$ is holomorphic on $\Lambda \cap \big( \mcal{B}_{\C}(u_*, r/4) \times \mcal{B}_{\C}((f_{a,*}, f_*), r/4) \big)$, we show that $h: z \mapsto \check{\mcal{D}}(u + z v, (f_a, f) + z(w_a, w))$
  is holomorphic on the open subset
  \[
    O_{h} = \{ z \in \C \ : \  (u, (f_a, f)) + z (v, (w_a, w)) \in \mcal{B}_{\C}(u_*, r/4) \times \mcal{B}_{\C}((f_{a,*}, f_*), r/4) \}.
  \]
  Let $z_0 \in O_{h}$, since $O_{h}$ is open there exists $\delta>0$ such that $C_{\delta} : [0,2\pi] \to O_{h}$ is a circle $\theta \mapsto z_0 + \delta e^{i\theta}$ of center $z_0$ and radius $\delta>0$ whose interior is contained in $O_{h}$.
  The function $\xi: [a,b] \to \C^n$ defined, using Cauchy's theorem on the holomorphic function $z\mapsto \check{D}(s, u(s) + z v(s), f(s) + z w(s))$, by
  \begin{align*}
    \xi(t)
     & = \int_{a}^{t}
    \frac{\partial \check{D}}{\partial x}(s, u(s)+z_0v(s), f(s)+z_0w(s)) v(s)\d s                                                  \\
     & \quad + \int_{a}^{t} \frac{\partial \check{D}}{\partial y}(s, u(s)+z_0v(s), f(s)+z_0w(s)) w(s)\d s                          \\
     & = \int_{a}^{t} \frac{1}{2\pi i} \oint_{C_{\delta}} \frac{\check{D}(s, u(s) + z v(s), f(s) + z w(s))}{(z - z_0)^2} \d z \d s
  \end{align*}
  is clearly absolutely continuous, i.e., $\xi \in \check{\mcal{U}}$.
  We now show that it is the complex derivative of $h$ at $z_0$ in $\check{\mcal{U}}$.
  To do this, we let $|\Delta z| < \delta$, and apply Cauchy's theorem to express the remainder of the difference quotient as
  \begin{align*}
    R_{\Delta z}(t) & \coloneqq
    \xi(t) - \frac{h(z_0+\Delta z)(t) - h(z_0)(t)}{\Delta z}                                \\
                    & = \int_{a}^{t} \frac{1}{2\pi i} \oint_{C_{\delta}}
    \check{D}(s, u(s) + z v(s), f(s) + z w(s))
    \left\{
    \frac{1}{(z - z_0)^2}\right.                                                            \\
                    & \hspace{5cm}\left.- \frac{1}{\Delta z}\left[
                                                              \frac{1}{z - (z_0 + \Delta z)}
                                                              - \frac{1}{z - z_0}
                                                              \right]
    \right\}
    \d z \d s.
  \end{align*}
  From this, we readily obtain
  \begin{align*}
     & \| R_{\Delta z} \|_{W^{1,1}([a,b];\C^n)}
    \leq                                                                                                                       \\
     & \left[\frac{1}{2}(b-a)^2 + (b-a) \right]\! M_{\mcal{C}} \delta \!\sup_{z \in C_{\delta}} \!\left| \frac{1}{(z - z_0)^2}
    - \frac{1}{\Delta z}\left[
                          \frac{1}{z - (z_0 + \Delta z)}
                          - \frac{1}{z - z_0}
                          \right] \right|
    \xrightarrow[\Delta z \to 0]{} 0,
  \end{align*}
  which proves that $\xi = \frac{d h}{d z}(z_0)$ is a G\^{a}teaux derivative of $\check{\mcal{D}}$ in $\AC([a,b];\C^n)$.
  Therefore, $h$ is holomorphic on $O_h$ and it follows that our extension $\check{\mcal{D}}$ is holomorphic on every affine slice $\Lambda \cap \big( \mcal{B}_{\C}(u_*, r/4) \times \mcal{B}_{\C}((f_{a,*}, f_*), r/4) \big)$.
  By \cite[Thm.~4]{bochnak1970analytic}, $\check{\mcal{D}}$ is Fr\'{e}chet holomorphic on the open neighborhood $\mcal{B}_{\C}(u_*, r/4) \times \mcal{B}_{\C}((f_{a,*}, f_*), r/4)$ of $(u_*, (f_{a,*}, f_*))$ in $\check{\mcal{U}} \times \check{\mcal{F}}$.
  Hence, the G\^{a}teaux derivative equals the Fr\'{e}chet derivative of $\check{\mcal{D}}$.
  The restriction of $\check{\mcal{D}}$ to the corresponding real neighborhood in $\mcal{U} \times \mcal{F}$ agrees with $\mcal{D}$, proving that it is real analytic and has the corresponding restricted Fr\'{e}chet derivative.
  The essentially uniform bound on the partial derivatives of $D$ is established using the same argument as in \cref{eqn:uniform_bound_on_derivatives_of_F}, which completes the proof.
\end{proof}

By linearity, \cref{lem:analytic_Nemytskii_operator} implies that the operator $\mcal{L}: \mcal{U} \times \mcal{F} \to \mcal{U}$ defined in \cref{eq:operator_L_dyn} is also analytic with analogous Fr\'{e}chet derivative and uniform bounds on $\partial_x F(s, u(s), f(s))$ and $\partial_y F(s, u(s), f(s))$.
We use the previous lemma together with the analytic implicit function theorem to prove that the solution operator provided by \cref{lem:continuity_of_ODE_solution} is analytic.
This enables us to apply the framework developed in \Cref{sec:analytic_systems} to analytic dynamical systems \eqref{eqn:analytic_ODE}.

\begin{lemma}[$\mathcal{G}$ is analytic]\label{lem:ODE_solution_mapping_is_analytic}
  Under \cref{ass:integrable_analytic_ODE}, the solution operator $\mcal{G}: \mcal{F} \to \mcal{U}$ defined in \cref{lem:continuity_of_ODE_solution} is analytic.
\end{lemma}
\begin{proof}
  By \cref{lem:analytic_Nemytskii_operator}, the operator
  $\mcal{A} = \mcal{D}_0 - \mcal{L}: \mcal{U} \times \mcal{F} \to \mcal{U}$ is analytic, and is given explicitly for $t\in[a,b]$ by
  $\mcal{A}[u, (f_a, f)](t) = u(t) - f_a - \int_{a}^{t} F(s, u(s), f(s)) \d s$.
  We observe that for all $(f_a, f) \in \mcal{F}$ the solution mapping $\mcal{G}:\mcal{F} \to \mcal{U}$ for \cref{eqn:analytic_ODE} satisfies $\mcal{A}[\mcal{G}(f_a, f), (f_a, f)] = 0$ .
  We prove that $\mcal{G}$ is analytic by applying the analytic implicit function theorem \cite{whittlesey1965analytic} to $\mcal{A}$.
  It suffices to show that at any point $(u_*, (f_{a,*}, f_*))$ with $u_* = \mcal{G}(f_{a,*}, f_*)$, the Fr\'{e}chet derivative of $\mcal{A}$ with respect to its first argument $\D_1 \mcal{A}[u_*, (f_{a,*}, f_*)] : \mcal{U} \to \mcal{U}$ has a bounded inverse. We do this by showing that the linear operator $T = I - \D_1 \mcal{A}[u_*, (f_{a,*}, f_*)] : \mcal{U} \to \mcal{U}$
  has operator norm smaller than $1$ with respect to an equivalent norm on $\mcal{U}$.
  By \cref{lem:analytic_Nemytskii_operator}, its action on $v \in \mcal{U}$ at $t \in [a,b]$ is given by
  \begin{align*}
    (T v)(t)
     & = \underbrace{v - \D_1 \mcal{D}_0 [u_*, (f_{a,*}, f_*)] v}_{0} + \D_1 \mcal{L} [u_*, (f_{a,*}, f_*)] v \\
    % = \D_1 \mcal{L} [u_*, (f_{a,*}, f_*)] v
     & = \int_{a}^{t} \frac{\partial F}{\partial x}(s, u_*(s), f_*(s)) v(s) \d s,
  \end{align*}
  and there exists a constant $L_u>0$ such that $\| \partial_x F(s, u_*(s), f_*(s))\| \leq L_u$ for almost every $s \in [a,b]$.
  Letting $\beta > 6 L_u$, we define an equivalent norm on $\mcal{U} = \AC([a,b];\R^n)$ by
  \[
    \vertiii{u}_{\beta}
    \coloneqq
    \sup_{t \in [a,b]} e^{-\beta (t-a)} \| u(t) \|
    + \int_{a}^{b} e^{-\beta (s-a)} \left\| \frac{\d u}{\d t}(s) \right\| \d s,
  \]
  which is equivalent to the $W^{1,1}([a,b];\R^n)$ norm by
  the Sobolev embedding from $W^{1,1}([a,b];\R^n) \hookrightarrow C([a,b];\R^n)$ and
  the fact that $0 < e^{-\beta (b-a)} \leq e^{-\beta (t-a)} \leq 1$ for all $t\in[a,b]$.
  \iffalse
    Hence, if $C$ is the Sobolev embedding constant, then
    \begin{multline*}
      \frac{e^{-\beta (b-a)}}{(1 + b - a)} \overbrace{
        \left(\| u \|_{L^1} + \left\| \frac{\d u}{\d t} \right\|_{L^1} \right)
      }^{\| u \|_{W^{1,1}}}
      % = \frac{1}{(1 + b - a)}\left[e^{-\beta (b-a)}\int_{a}^{b} \| u(s) \| \d s + e^{-\beta (b-a)}\int_{a}^{b} \| \frac{\d u}{\d t}(s) \| \d s\right] \\
      \leq\\ \frac{1}{(1 + b - a)}\left[(b-a) e^{-\beta (b-a)} \sup_{t \in [a,b]} \| u(t) \| + e^{-\beta (b-a)} \int_{a}^{b} \left\| \frac{\d u}{\d t}(s) \right\| \d s\right] \\
      \leq \underbrace{
        \sup_{t \in [a,b]} e^{-\beta (t-a)} \| u(t) \|
        + \int_{a}^{b} e^{-\beta (s-a)} \left\| \frac{\d u}{\d t}(s) \right\| \d s
      }_{\vertiii{u}_{\beta}}
      \leq \sup_{t \in [a,b]} \| u(t) \|
      + \int_{a}^{b} \left\| \frac{\d u}{\d t}(s) \right\| \d s\\
      \leq (C + 1) \underbrace{
        \left(\| u \|_{L^1} + \left\| \frac{\d u}{\d t} \right\|_{L^1} \right)
      }_{\| u \|_{W^{1,1}}}.
    \end{multline*}
  \fi
  Choosing $v \in \mcal{U}$, combining the expression of $T$ and the definition of the equivalent norm yields
  \[
    \vertiii{T v}_{\beta}
    \leq L_u \sup_{t\in [a,b]} e^{-\beta (t-a)} \int_{a}^{t} \| v(s) \| \d s
    + L_u \int_{a}^{b} e^{-\beta (t-a)} \left\| v(t) \right\| \d t.
  \]
  We control the first term using the bound $\| v(s) \| \leq e^{\beta (s-a)} \vertiii{v}_{\beta}$, and the second term with $\| v(t) \| \leq \| v(a) \| + \int_{a}^{t} \left\| \frac{\d v}{\d t}(s) \right\| \d s$ (derived from the fundamental theorem of calculus) where $\| v(a) \| \leq \vertiii{v}_{\beta}$, to obtain
  \begin{align*}
    \vertiii{T v}_{\beta}
     & \leq L_u \vertiii{v}_{\beta} \sup_{t\in [a,b]} e^{-\beta (t-a)} \underbrace{
                                                                         \int_{a}^{t} e^{\beta (s-a)} \d s
                                                                       }_{\leq  e^{\beta (t - a)}/\beta} \\&\quad
    + L_u \vertiii{v}_{\beta}
    \underbrace{
      \int_{a}^{b} e^{-\beta (t-a)} \d t
    }_{\leq 1/\beta}
    + L_u \int_{a}^{b} \int_{a}^{t} e^{-\beta (t - a)} \left\| \frac{\d v}{\d t}(s) \right\| \d s \d t    \\
     & \leq \frac{2 L_u}{\beta} \vertiii{v}_{\beta}
    + L_u \int_{a}^{b}
    \underbrace{
      \int_{s}^{b} e^{-\beta (t - a)} \d t
    }_{\leq  e^{-\beta (s-a)}/\beta}
    \left\| \frac{\d v}{\d t}(s) \right\| \d s
    \leq \frac{3 L_u}{\beta} \vertiii{v}_{\beta} < \frac{1}{2}\vertiii{v}_{\beta},
  \end{align*}
  where we used Fubini's theorem in the last integral term of the first inequality. Therefore, the operator norm $\vertiii{T}_{\beta} = \sup_{v \in \mcal{U}, v \neq 0} \vertiii{T v}_{\beta} / \vertiii{v}_{\beta} < 1/2$.
  It follows that $\D_1 \mcal{A}(u_*, (f_{a,*}, f_*)) = I - T$ has a bounded inverse, and so $\mcal{G}$ is analytic by the implicit function theorem~\cite{whittlesey1965analytic}.
\end{proof}

Now that we have verified the required analyticity of the dictionary functions and analytic well-posedness of the solution operator, we apply the framework from \Cref{sec:analytic_systems} to characterize the inputs $(f_a, f) \in \mcal{F}$ allowing the coefficients $\vect{c}^*$ of the underlying dynamical system \eqref{eqn:analytic_ODE} to be identified. We consider the linear operator $L_{(f_a, f), u} : \R^N \to \mcal{U}$, defined similarly to  \cref{eqn:analytic_system_recoverability_operator} using the integrated dictionary functions in \cref{eqn:analytic_ODE_dict_ops} as
\begin{equation}\label{eqn:analytic_ODE_recoverability_operator}
  \left(L_{(f_a, f), u} \vect{c}\right)(t) =  \sum_{i = 1}^{N} c_i \int_{a}^{t} D_i(s, u(s), f(s)) \d s,\quad \vect{c}\in\R^N,\quad t\in[a,b].
\end{equation}
Given the input $(f_a,f)\in\mathcal{F}$ and solution $u\in\mathcal{U}$, the unknown coefficients in \cref{eqn:analytic_ODE} are the unique solution $\vect{c} = \vect{c}^*$ to the linear system $
  L_{(f_a, f), u} (\vect{c}) = u - f_a$
when $L_{(f_a, f), u}$ is injective. Here, we recover the method proposed by \cite{schaeffer2017sparse} to identify the unknown coefficients in the integral form of \cref{eqn:analytic_ODE} given by \cref{eq:operator_L_dyn}. The degenerate set of inputs for which coefficients cannot be recovered is given by
\begin{equation}\label{eqn:degenerate_set_for_analytic_ODE}
  \Fdeg(\vect{c}^*) = \left\{ (f_a, f) \in \mcal{F} \mid \ L_{(f_a, f), \mcal{G}(f_a, f)} \ \mbox{not injective} \right\}.
\end{equation}
The following corollaries are a consequence of \cref{cor:analytic_recovery,cor:analytic_recovery_alternative} and provide a characterization of the degenerate set for analytic dynamical systems.

\begin{corollary}[Degenerate set of analytic dynamical systems]\label{cor:analytic_ODE_recovery}
  Suppose that $\mcal{F}$ has a countable Schauder basis and \cref{ass:integrable_analytic_ODE} holds.
  If there exists $(f_{a,0}, f_0) \in \mcal{F} \setminus \Fdeg(\vect{c}^*)$, then
  $\Fdeg(\vect{c}^*)$ is compactly slice-null and has measure zero with respect to any nondegenerate Gaussian measure on Borel sets of $\mcal{F}$.
\end{corollary}

\begin{corollary}[Zero-one law for analytic dynamical systems]\label{cor:analytic_recovery_alternative_for_ODE}
  Suppose that $\mcal{F}$ has a countable Schauder basis and \cref{ass:integrable_analytic_ODE} holds. Then, exactly one of the following is true:
  \begin{enumerate}[label=\roman*), itemsep=0em, leftmargin=*]
    \item If $(f_a,f)$ is sampled from a nondegenerate Gaussian measure on $\mcal{F}$, then $(f_a,f) \notin \Fdeg(\vect{c}^*)$ with probability $1$. \label{recov_1_dynamical}
    \item $\Fdeg(\vect{c}^*) = \mcal{F}$.
  \end{enumerate}
\end{corollary}

\begin{remark}[Existence of Schauder basis]
  \Cref{cor:analytic_ODE_recovery,cor:analytic_recovery_alternative_for_ODE} assume that $\mcal{F}$ has a countable Schauder basis, which is not true for $L^{\infty}([a,b];\R^m)$.
  However, countable Schauder bases do exist for certain closed subspaces of $L^{\infty}([a,b];\R^m)$, such as the space of continuous functions $C([a,b];\R^m)$, or bounded piecewise continuous functions with finitely many fixed intervals of continuity.
\end{remark}

The following example illustrates both possible outcomes of the recovery alternative in \cref{cor:analytic_recovery_alternative_for_ODE}.
\begin{example}[Uncontrollable linear system]    \leavevmode
  \begin{enumerate}[label=\roman*), itemsep=0em, leftmargin=*]
    \item Consider the linear dynamical system
          \[
            \frac{\d}{\d t}
            \begin{bmatrix}
              u_1(t) \\
              u_2(t)
            \end{bmatrix}
            =
            \begin{bmatrix}
              -1 & c_1^* \\
              0  & c_2^*
            \end{bmatrix}
            \begin{bmatrix}
              u_1(t) \\
              u_2(t)
            \end{bmatrix}
            +
            \begin{bmatrix}
              1 \\
              0
            \end{bmatrix}
            f(t),
            \qquad
            \begin{bmatrix}
              u_1(0) \\
              u_2(0)
            \end{bmatrix}
            =
            \begin{bmatrix}
              0 \\
              0
            \end{bmatrix},
          \]
          with unknown coefficients $\vect{c}^* = (c_1^*, c_2^*)$. If the inputs are selected from the subspace $\mcal{F} = \left\{ (0, f) : f \in C([0,1]; \R) \right\} \subset \R^2\times C([0,1];\R)$ with zero initial conditions, then the unknown coefficients cannot be determined because the state variable $u_2(t)$ remains zero for all time $t$ irrespective of the choice of $f$, meaning $\Fdeg = \mcal{F}$.
    \item On the other hand, if we consider the autonomous dynamical system determined by its initial conditions:
          \[
            \frac{\d}{\d t}
            \begin{bmatrix}
              u_1(t) \\
              u_2(t)
            \end{bmatrix}
            =
            \begin{bmatrix}
              -1 & c_1^* \\
              0  & c_2^*
            \end{bmatrix}
            \begin{bmatrix}
              u_1(t) \\
              u_2(t)
            \end{bmatrix}
            ,
            \qquad
            \begin{bmatrix}
              u_1(0) \\
              u_2(0)
            \end{bmatrix}
            =
            \begin{bmatrix}
              f_{0,1} \\
              f_{0,2}
            \end{bmatrix},
          \]
          with input space $\mcal{F} = \left\{ (f_0,0) : f_0\in \R^2\right\} \subset \R^2\times L^{\infty}([0,1]; \R)$, then
          every input with $f_{0,2} \neq 0$
          produces a solution $u$ that can be used to identify $c_1^*$ and $c_2^*$. \cref{cor:analytic_ODE_recovery} ensure that $\Fdeg$ is compactly slice-null and has zero nondegenerate Gaussian measure in $\mcal{F}$.
  \end{enumerate}
\end{example}

\section{Recovery with Differential Operators}
\label{sec:differential_operators}
The recovery result for analytic systems stated in \cref{cor:analytic_recovery} assumes the existence of a recovering input $f_0 \in \mcal{F} \setminus \Fdeg(\vect{c}^*)$. This section explains how certain analytic systems of differential equations are guaranteed to have a recovering input, meaning the conclusions of \cref{cor:analytic_recovery} automatically hold. We consider a particular case of \cref{eqn:analytic_system} with systems of the form
\begin{equation}\label{eqn:Lu_f_differential_ops}
  \mcal{L}(u)\coloneqq \sum_{i=1}^{N} c^*_i \mcal{D}_i(u) = f,
\end{equation}
where $\mcal{L},\ \mcal{D}_i: \mcal{U} \to \mcal{F}$ are analytic operators between Banach spaces $\mcal{U}$ and $\mcal{F}$ of $\R^n$-valued distributions on an open set $\Omega \subset \R^d$.
Specifically, each $f \in \mcal{F}$ is an $n$-tuple of distributions $f = (f_1, \ldots, f_n)$ whose action on a test function $\phi \in C^{\infty}_c(\Omega)$ is defined component-wise as $f[\phi] = (f_1[\phi], \ldots, f_n[\phi]) \in \R^n$. Examples of Banach spaces of $\R^n$-valued distributions include $L^p(\Omega;\R^n)$, for $1 \leq p \leq \infty$, and negative Sobolev spaces such as $H^{-s}(\Omega;\R^n)$, for $s > 0$.
Clearly \cref{eqn:Lu_f_differential_ops} is a special case of \cref{eqn:analytic_system}, and we make the following analytic well-posedness assumption about the underlying system we aim to identify:
\begin{assumption}[Well-posedness of differential equations]\label{ass:diff_op_well_posed}
  We assume that \eqref{eqn:Lu_f_differential_ops} is well-posed for $f\in\mcal{F}$ with analytic solution operator $\mcal{G} = \mcal{L}^{-1}: \mcal{F} \to \mcal{U}$.
\end{assumption}

For $u \in C^{\infty}_c(\Omega;\R^n)$ and a multi-index $\alpha = (\alpha_1, \ldots, \alpha_d) \in \mathbb{N}_0^{d}$, we define the partial derivative $\partial^{\alpha} u = \partial_{x_1}^{\alpha_1} \cdots \partial_{x_d}^{\alpha_d} u$, where $|\alpha| = \alpha_1 + \cdots + \alpha_d$.
The \emph{$r$-jet} of $u$ is the function $\jet^r u : \Omega \to \mcal{J}(r,d,n) \coloneqq \R^{n\binom{d+r}{d}}$ consisting of $u$ and all of its partial derivatives up to order $r$: $\jet^r u(x) = \big[ \partial^{\alpha} u(x) \big]_{|\alpha|\leq r}$ for $x \in \Omega$. We assume that the dictionary operators $\mcal{D}_i$ are given by integrating functions of the $r$-jet of $u$ against test functions, which is a common form for dictionary operators in systems of differential equations~\cite{messenger2021weak}.

\begin{assumption}[Choice of dictionary operators]\label{ass:analytic_diff_op}
  Suppose $C^{\infty}_c(\Omega;\R^n) \subset \mcal{U}$ and for $1\leq i\leq N$, $\mcal{D}_i:\mcal{U} \to \mcal{F}$ is analytic as an operator between Banach spaces, and that
  there is an integer $r$ and continuous functions $D_i : \Omega \times \mcal{J}(r,d,n) \to \R^n$ such that
  \begin{equation}\label{eqn:diff_dict_funs}
    \mcal{D}_i(u)[\phi] = \int_{\Omega} D_i(x, \jet^r u(x)) \phi(x) \d x,\quad u \in C^{\infty}_c(\Omega;\R^n),\quad \phi \in C^{\infty}_c(\Omega).
  \end{equation}
\end{assumption}

For example, \cref{ass:analytic_diff_op,ass:diff_op_well_posed} hold for the Poisson problem $\mcal{L}(u) = u_{xx} + u_{yy} = f$ on the unit square $\Omega = (0,1)^2 \subset \R^2$ with spaces $\mcal{U} = H^2(\Omega) \cap H_0^1(\Omega)$ and $\mcal{F} = L^2(\Omega)$, and dictionary operators $u \mapsto u, u^2, u^3, \sin(u), e^u, u_{x}, u_{y}, u_{xx}, u_{xy}, u_{yy}$.
The nonlinear dictionary operators are analytic thanks to \cite[Ex.~4.3.6]{buffoni2003analytic} together with the Sobolev embedding $H^2(\Omega) \hookrightarrow C(\Omega)$ and continuous linear embedding $C(\Omega) \hookrightarrow L^2(\Omega)$.

We show that a recovering input exists by constructing a smooth solution $u_0 \in \mcal{U}$ such that the dictionary operators take linearly-independent values across a collection of distinct locations in the spatial domain $\Omega$.
This yields a recovery algorithm for the unknown coefficients that works without \cref{ass:diff_op_well_posed} in the case where $u \in \mcal{U}$ can be selected and we observe $\mcal{L}(u)$.
A recovery algorithm for the case where $f \in \mcal{F}$ can be selected and we observe $u = \mcal{G}(f)$ is then obtained from \cref{cor:analytic_recovery} and the existence of a recovering input $f_0 = \mcal{L}(u_0)$ found by applying the differential operator.

The following technical lemma is used in our construction of $u_0$.
Specifically, applying the lemma to $D_1, \ldots, D_N$ allows us to find a set of spatial locations in $\Omega$ and $r$-jets, and construct a function $u_0$ by interpolating at these locations to ensure that $\mcal{D}_1(u_0), \ldots, \mcal{D}_N(u_0)$ are linearly independent, such that the unknown coefficients can be recovered using the general method from \Cref{sec:analytic_systems}.

\begin{lemma}[Linear independence and injectivity]\label{lem:existence_of_sample_points}
  Let $m\geq 1$, then functions $f_1, \ldots, f_m : \mcal{X} \to \mcal{Y}$ from a set $\mcal{X}$ into a real vector space $\mcal{Y}$ are linearly independent if and only if there exist $x_1, \ldots, x_m \in \mcal{X}$ such that the quasimatrix
  \[
    A_m(x_1, \ldots, x_m) =
    \begin{bmatrix}
      f_1(x_1) & \cdots & f_m(x_1) \\
      \vdots   & \ddots & \vdots   \\
      f_1(x_m) & \cdots & f_m(x_m)
    \end{bmatrix}
    : \R^m \to \mcal{Y}^m
  \]
  is injective.
\end{lemma}
\begin{proof}
  We first suppose that there exist $x_1, \ldots, x_m$ such that $A_m(x_1, \ldots, x_m)$ is injective.
  Let $\vect{c} = (c_1, \ldots, c_m) \in \R^m$ be such that $\sum_{i=1}^{m} c_i f_i = 0$. This implies that $A_m(x_1, \ldots, x_m) c = 0$, which by injectivity of $A_m(x_1, \ldots, x_m)$ implies that $c = 0$, and $f_1, \ldots, f_m$ are linearly independent.

  We now assume that the family of functions is linearly independent and proceed by induction over $m$.
  If $m=1$, since $f_1$ is nonzero, there exists $x_1 \in \mcal{X}$ such that $A_1(x_1) = [f_1(x_1)] \neq 0$.
  Suppose that there exists $x_1, \ldots, x_{m} \in \mcal{X}$ such that $A_{m}(x_1, \ldots, x_{m})$ is injective. Proceeding by contradiction, suppose that $A_{m+1}(x_1, \ldots, x_{m}, x)$ fails to be injective for every $x \in \mcal{X}$. Then, for every $x \in \mcal{X}$, there exists a nonzero vector $\vect{c}(x) = (c_1(x), \ldots, c_{m+1}(x)) \in \R^{m+1}$ such that $A_{m+1}(x_1, \ldots, x_{m}, x) \vect{c}(x) = 0$. The first $m$ rows of this equation read
  \[
    A_{m}(x_1, \ldots, x_{m})
    \begin{bmatrix}
      c_1(x) \\
      \vdots \\
      c_m(x)
    \end{bmatrix} =- c_{m+1}(x)
    \begin{bmatrix}
      f_{m+1}(x_1) \\
      \vdots       \\
      f_{m+1}(x_m)
    \end{bmatrix}.
  \]
  Since $A_{m}(x_1, \ldots, x_{m})$ is injective, $c_{m+1}(x)$ uniquely determines $c_1(x),\ldots c_m(x)$, and $\Null A_{m+1}(x_1, \ldots, x_{m}, x)$ is $1$-dimensional.
  Moreover, it is contained in the nullspace of $:A_{m+1}(x_1, \ldots, x_{m}): \R^{m+1} \to \mcal{Y}^m$. Hence, $\Null A_{m+1}(x_1, \ldots, x_{m}, x)$ does not depend on $x$. Letting $\vect{c} = (c_1, \ldots, c_{m+1}) \in \R^{m+1}$ span this nullspace, we have $A_{m+1}(x_1, \ldots, x_{m}, x) \vect{c} = 0$ for every $x \in \mcal{X}$, whose last row reads
  \[
    c_1 f_1(x) + \cdots + c_m f_m(x) + c_{m+1} f_{m+1}(x) = 0, \qquad \forall x \in \mcal{X},
  \]
  contradicting the linear independence of these functions.
  Therefore, there exists $x_{m+1}\in\mcal{X}$ such that $A_{m+1}(x_1, \ldots, x_{m}, x_{m+1})$ injective.
\end{proof}

Applying \cref{lem:existence_of_sample_points} to the functions $D_1, \ldots, D_N$ from \cref{eqn:diff_dict_funs} provides a collection of spatial locations and $r$-jets, which if assigned to $u_0$ guarantees linear independence of $\mcal{D}_1(u_0), \ldots, \mcal{D}_N(u_0)$.
\Cref{prop:analytic_diff_op_u0} makes this rigorous by constructing $u_0$ with the prescribed $r$-jets, and showing that it enables coefficient recovery by making the operator in \cref{eqn:analytic_system_recoverability_operator} injective.

\begin{proposition}[Existence of $u_0$]\label{prop:analytic_diff_op_u0}
  Suppose \cref{ass:analytic_diff_op} holds with linearly independent functions $\{ D_i \}_{i=1}^{N}$, then there exists $u_0 \in C_c^{\infty}(\Omega;\R^n) \subset \mcal{U}$ such that the linear operator $L_{u_0} : \R^N \to \mcal{F}$ defined by $
    L_{u_0}(\vect{c}) =  \sum_{i=1}^{N} c_i \mcal{D}_i(u_0)$ is injective.
\end{proposition}

\begin{proof}
  To construct $u_0$, we first observe that by \cref{lem:existence_of_sample_points}
  there exist points $\{(\bar{x}_i, \vect{\xi}_i)\}_{i=1}^{N} \subset \Omega \times \mcal{J}(r,d,n)$ such that the matrix
  $
    [D_j(\bar{x}_i, \vect{\xi}_i)]_{1 \leq i,j \leq N} \in \R^{n N \times N}
  $
  is injective. Since the functions $D_j$ are continuous, the elements $\{\bar{x}_i\}_{i=1}^{N} \subset \Omega$ can be perturbed so that they are all distinct while maintaining injectivity of $[D_j(\bar{x}_i, \vect{\xi}_i)]_{1 \leq i,j \leq N}$.
  Let $\{ B(\bar{x}_i, r_i) \}_{i=1}^{N}$ be open balls centered at each $\bar{x}_i$ with radii $r_i > 0$ sufficiently small such that every $B(\bar{x}_i, r_i) \subset \Omega$ and the balls are disjoint, i.e., $B(\bar{x}_i, r_i) \cap B(\bar{x}_j, r_j) = \emptyset$ when $i \neq j$.

  For $1\leq i\leq N$, we write $\vect{\xi}_i =(\xi_\alpha)_{|\alpha|\leq r}$ and consider the $\R^n$-valued polynomial $\varphi_i(x) = \sum_{|\alpha|\leq r}\tfrac{\xi_\alpha}{\alpha!}(x-\bar{x}_i)^\alpha$ on $\R^d$ whose $r$-jets satisfy $\jet^r \varphi_i(\bar{x}_i) = \vect{\xi}_i$. Let $\theta_i \in C_c^{\infty}(\R^d)$ be bump functions satisfying $\theta_i(x) = 1$ for all $x\in B(\bar{x}_i, r_i/2)$ with compact support such that $\supp(\theta_i) \subset B(\bar{x}_i, r_i)$.
  We observe that the function
  $
    u_0 = \sum_{i=1}^N \theta_i \varphi_i\in C_c^{\infty}(\Omega;\R^n) \subset \mcal{U}$ and
  $
    \jet^r u_0(\bar{x}_i) = \jet^r \varphi_i(\bar{x}_i) = \vect{\xi}_i$ for $1\leq i\leq N$.
  Therefore, the $(n N) \times N$ matrix
  \[
    \mat{\bar{A}} = [D_j(\bar{x}_i, \jet^r u_0(\bar{x}_i))]_{1 \leq i,j\leq N}
    = [D_j(\bar{x}_i, \vect{\xi}_i)]_{1 \leq i,j\leq N}
  \]
  is injective.

  The rank function is lower semicontinuous, hence the injective real $(n N) \times N$ matrices are an open subset of $\R^{(n N) \times N}$. Therefore, there exists $\varepsilon > 0$ such that if $\mat{A} \in \R^{(n N) \times N}$ satisfies $\| \mat{A} - \mat{\bar{A}} \|_{\infty}  < \varepsilon$, then $\mat{A}$ is injective. Moreover, since each function $x \mapsto D_j(x, \jet^r u_0(x))$ is continuous, there exists $\delta > 0$ such that for every $1 \leq i,j \leq N$, if $\| x - \bar{x}_i \| < \delta$, then $\| D_j(x, \jet^r u_0(x)) - D_j(\bar{x}_i, \jet^r u_0(\bar{x}_i)) \| < \varepsilon$. Let $\phi_i \in C_c^{\infty}(\Omega)$ be nonnegative with support contained in $B(\bar{x}_i, \delta)$, normalized such that $\int_{\Omega} \phi_i(x) \d x = 1$ for $1\leq i \leq N$.
  Then, we have
  \begin{multline*}
    \left\| \int_{\Omega} D_j(x, \jet^r u_0(x)) \phi_i(x) \d x - D_j(\bar{x}_i, \jet^r u_0(\bar{x}_i)) \right\|\\
    \leq
    \int_{\Omega} \left\| D_j(x, \jet^r u_0(x)) - D_j(\bar{x}_i, \jet^r u_0(\bar{x}_i)) \right\| \phi_i(x) \d x
    < \varepsilon,
  \end{multline*}
  meaning $n N \times N$ matrix defined by $\mat{A} = \left[ \int_{\Omega} D_j(x, \jet^r u_0(x)) \phi_i(x) \d x \right]_{1 \leq i, j \leq N}$ is injective.
  Since $x \mapsto D_j(x, \jet^r u_0(x))$ agrees with $\mcal{D}_j(u_0)$ as an $\R^n$-valued distribution by \cref{ass:analytic_diff_op}, it follows that the operator $\R^N \to \R^{n N}$ given by
  \[\resizebox{\textwidth}{!}{$
  \vect{c}
  \longmapsto
  \begin{bmatrix}
    (L_{u_0} \vect{c})[\phi_1] \\
    \vdots                     \\
    (L_{u_0} \vect{c})[\phi_N]
  \end{bmatrix}
  =
  \begin{bmatrix}
    \int_{\Omega} D_1(x, \jet^r u_0(x)) \phi_1(x) \d x
           & \cdots &
    \int_{\Omega} D_N(x, \jet^r u_0(x)) \phi_1(x) \d x \\
    \vdots & \ddots & \vdots                           \\
    \int_{\Omega} D_1(x, \jet^r u_0(x)) \phi_N(x) \d x
           & \cdots &
    \int_{\Omega} D_N(x, \jet^r u_0(x)) \phi_N(x) \d x
  \end{bmatrix}
  \begin{bmatrix}
    c(1)   \\
    \vdots \\
    c(N)
  \end{bmatrix}
  = \mat{A} \vect{c}$}
  \]
  is injective. Finally, suppose $L_{u_0}\vect{c}=0$; evaluating at $(\phi_1,\ldots,\phi_N)$ gives $\mat{A} \vect{c}=0$, hence $\vect{c}=0$ by injectivity of $\vect{c}\mapsto\mat{A} \vect{c}$. Therefore $L_{u_0}$ is injective, which completes the proof.
\end{proof}

The following corollary to \cref{prop:analytic_diff_op_u0} invokes results from \Cref{sec:analytic_systems} to establish that the set of $u\in \mcal{U}$ failing to make $L_u$ injective is compactly slice null and has zero nondegenerate Gaussian measure.
This provides a recovery algorithm for the case in which one can select $u$ and observe $f = \mcal{L}(u)$. That is, selecting $u$ from a nondegenerate Gaussian measure yields $\vect{c} = \vect{c}^*$ as the unique solution of $L_u \vect{c} = \mcal{L}(u)$ with probability $1$. To use this algorithm, one does not have to assume that $\mcal{L}^{-1}$ exists. We first introduce the degenerate solution set.

\begin{definition}[Degenerate solution set]\label{def:degenerate_sol_set}
  We define the degenerate solution set associated with \cref{eqn:Lu_f_differential_ops} as $
    \Udeg =
    \left\{ u \in \mcal{U} \mid \exists\vect{a}\in\R^N\setminus\{0\},\, L_{u}(\vect{a}) = 0\right\}$.
\end{definition}

Note that unlike the degenerate input set $\Fdeg(\vect{c}^*)$, the degenerate solution set $\Udeg$ does not depend on the unknown coefficients $\vect{c}^*$, and characterizes the recovery of all models that are linearly parameterized by the dictionary $\{ \mcal{D}_i \}_{i=1}^{N}$.

\begin{corollary}[Prevalence of degenerate solutions]\label{cor:analytic_diff_op_u_degenerate}
  Suppose $\mcal{U}$ has a countable Schauder basis and \cref{ass:analytic_diff_op} holds with linearly independent functions $\{ D_i \}_{i=1}^{N}$, then $\Udeg$ is compactly slice-null and has measure zero with respect to any nondegenerate Gaussian measure on Borel sets of $\mcal{U}$.
\end{corollary}
\begin{proof}
  The proof consists of applying \cref{cor:analytic_recovery} after swapping the spaces $\mcal{F}$ and $\mcal{U}$. Specifically, we consider the system $\mcal{L}(\tilde{f}) = \tilde{u}$ with $\tilde{f} \in \tilde{\mcal{F}} \coloneqq \mcal{U}$ and $\tilde{u} \in \tilde{\mcal{U}} \coloneqq \mcal{F}$.
  By \cref{ass:analytic_diff_op}, $\tilde{\mcal{G}} = \mcal{L}$ is the analytic solution operator, and following \cref{prop:analytic_diff_op_u0} there exists $\tilde{f}_0 \equiv u_0 \in \mcal{U} \setminus \Udeg = \mcal{\tilde{F}}\setminus \mcal{\tilde{F}}_{\rm{deg}}(\vect{c}^*)$. The conclusion follows from \cref{cor:analytic_recovery}.
\end{proof}

We now explain why coefficient recovery is impossible for certain spaces $\mcal{U}$ and dictionaries that violate \cref{ass:analytic_diff_op}, showing the importance of the assumptions. Hence, the assumption $C_c^\infty(\Omega)\subset \mcal U$ rules out such globally constrained classes and is precisely what allows one to prescribe arbitrary jets at finitely many points.

\begin{example}[Importance of \cref{ass:analytic_diff_op}]
  There are some situations where $\Udeg=\mcal U$ can occur, i.e., coefficient recovery is impossible, such as when the admissible class $\mcal U$
  is itself constrained by a differential identity lying in the span of the dictionary.
  \begin{itemize}[leftmargin=*]
    \item A first example is the harmonic class on $\Omega\subset\R^2$,
          \[
            \mcal U=\{u\in C^2(\Omega)\mid \Delta u=0\},
            \qquad \text{with} \quad
            \mcal D_1(u)=u_{xx},\quad \mcal D_2(u)=u_{yy},
          \]
          In this case, $\mcal D_1(u)+\mcal D_2(u)=0$ for all $u\in\mcal U$, therefore $L_u(\vect{c}) = (c_1-c_2)u_{xx}$ and one cannot distinguish the coefficients $\vect{c}\in \vspan\{(1,1)\}$. Hence, every $u\in\mcal U$ is degenerate and $\Udeg=\mcal U$. Note that the set defined as $\Udeg^2 = \left\{ u \in \mcal{U} \mid \rank L_u =0 \right\}$ is given locally on a rectangular domain $\Omega$ by $\Udeg^2=\vspan\{1,x,y,xy\}$. Here, one can recover the coefficients up to a one-dimensional subspace for $u\notin \Udeg^2$.
    \item A second example is the fixed-frequency trigonometric space on an interval $I\subset \R$,
          \[
            \mcal U=\vspan\{x\mapsto \sin(kx), x\mapsto \cos(kx)\},
            \qquad \text{with} \quad
            \mcal D_1(u)=u,\quad \mcal D_2(u)=u_{xx},
          \]
          for which $\mcal D_2(u)+k^2 \mcal D_1(u)=0$ for all $u\in\mcal U$. Then, $L_u(\vect{c}) = c_1 u + c_2 u_{xx} = (c_1-k^2c_2) u$ for all $u\in\mcal U$, and $\Udeg=\mathcal{U}$, $\Udeg^2=\{0\}$ so one can always recover the coefficients up to a one-dimensional subspace for $u\neq 0$.
    \item A third example is the space of affine functions on $\Omega\subset\R^d$. If the dictionary contains any second-order derivative, then that dictionary element vanishes identically on $\mcal U$, so every $u\in\mcal U$ is degenerate.
  \end{itemize}
\end{example}

The following corollary of \cref{prop:analytic_diff_op_u0} is a coefficient recovery guarantee for \cref{eqn:Lu_f_differential_ops} when $f \in \mcal{F}$ can be chosen and we observe the resulting solution $u = \mcal{G}(f)$.
It is only at this point that we actually require the analytic well-posedness from \cref{ass:diff_op_well_posed}.
The resulting recovery algorithm proceeds as follows:
when $f$ is sampled from a nondegenerate Gaussian measure and $u = \mcal{G}(f)$, the correct coefficients $\vect{c} = \vect{c}^*$ will be the unique solution of $L_{u} \vect{c} = f$ with probability $1$.

\begin{corollary}[Almost sure recovery of differential operators] \label{cor:analytic_recovery_strong}
  Suppose that $\mcal{F}$ has a countable Schauder basis, \cref{ass:diff_op_well_posed} holds, and \cref{ass:analytic_diff_op} holds with linearly independent functions $\{ D_i \}_{i=1}^{N}$.
  Then,
  \[
    \Fdeg(\vect{c}^*)=
    \left\{ f \in \mcal{F} \mid \exists\vect{a}\in\R^N\setminus\{0\},\, L_{\mcal{G}(f)}(\vect{a}) = 0 \right\}
    = \mcal{L}(\Udeg),
  \]
  is compactly slice-null and has measure zero with respect to any nondegenerate Gaussian measure on Borel sets of $\mcal{F}$.
\end{corollary}

\begin{proof}
  By \cref{prop:analytic_diff_op_u0} there exists $u_0 \in \mcal{U}$ such that $L_{u_0}$ is injective. Taking $f_0 = \mcal{L}(u_0)$ makes $L_{\mcal{G}(f_0)}$ injective. Therefore, \cref{cor:analytic_recovery} yields the conclusion.
\end{proof}

Note that \cref{cor:analytic_recovery_strong} is stronger than \cref{prop:analytic_diff_op_u0} as we guarantee the existence of $f_0\in\mcal{F}$ that allows recovery of the differential operator. In this case, we show almost sure recovery of the coefficients when $f$ is sampled from a nondegenerate Gaussian measure on $\mcal{F}$ and $u = \mcal{G}(f)$ is observed. \Cref{sec:second_order_PDE} provides an example where the degenerate set is non-trivial.

\section{Refined characterizations of degenerate sets}
\label{sec:refined_characterization}

\cref{thm:degenerate_sets} ensures that the degenerate input set is either compactly slice null or equals the entire input set. In this section, we provide a refined characterization of the degenerate sets for single-variable polynomial ODEs and specific partial differential equations.

\subsection{Single-variable polynomial ODEs}
We use a method that differs from the framework introduced in \Cref{sec:analytic_systems} to provide detailed information about the degenerate sets for polynomial ordinary differential equations with a single dependent variable.
For these systems we show that functions in the degenerate sets have guaranteed regularity properties, and thus constructing inputs lacking this regularity allows to recover the unknown coefficients.

Let $I \subset \R$ be a compact interval with nonempty interior, and consider an ordinary differential equation on $I$ of the form
\begin{equation}\label{eqn:nonlinear_ODE}
  \mcal{L}(u)(t) = \sum_{i=1}^{N_1}\sum_{j=1}^{N_2} c^*_{i,j} \underbrace{\theta_{i}(t) \psi_{j}( \jet^r u(t))}_{\mcal{D}_{i,j}(u)(t)} = f(t), \quad t \in I,
\end{equation}
where $f \in \mcal{F} = C(I)$ is a known continuous scalar-valued input function, $\vect{c}^* \in \R^{N_1\times N_2}$ is a vector of unknown coefficients, $\theta_{i} \in C^{\infty}(I)$ are known linearly-independent smooth functions, and $\psi_{j} : \mcal{J}(r,1,1) \to \R$ are known linearly-independent polynomials on the jet space $\mcal{J}(r,1,1) = \R^{r+1}$. We consider only strong solutions of \cref{eqn:nonlinear_ODE}, satisfying $u \in \mcal{U} = C^{r}(I)$, and assume that such a solution exists for each $f \in C(I)$.
Solutions need not be unique since we leave the boundary conditions unspecified.
The results we develop apply to any solution regardless of boundary conditions.

We aim to recover $\vect{c}^*$ using a single input-solution pair $(f,u) \in C(I) \times C^{r}(I)$ and define the parameter-to-output recovery map $L_{u} : \R^{N_1\times N_2} \to C(I)$ as
\[
  L_u : \vect{c} \longmapsto \sum_{i=1}^{N_1}\sum_{j=1}^{N_2}  c_{i,j} \mcal{D}_{i,j}(u).
\]
If $L_u$ is injective, then solving the linear system $L_u \vect{c} = f$ recovers the correct coefficients $\vect{c} = \vect{c}^*$.
This fails whenever $u$ lies in the degenerate solution set
\[
  \Udeg \coloneqq
  \left\{ u \in C^r(I) \mid \mbox{$L_u$ is not injective} \right\}.
\]
The inputs $f$ that can produce a degenerate solution lie in the set
\[
  \Fdeg(\vect{c}^*)
  \coloneqq \left\{ f \in C(I) \mid \mbox{ $\exists u \in \mcal{L}^{-1}(f)$ s.t. $L_u$ is not injective} \right\}
  = \mcal{L}(\Udeg).
\]
We now prove a regularity theorem for functions in $\Udeg$ and $\Fdeg(\vect{c}^*)$.

\begin{theorem}[Regularity of degenerate solutions]\label{thm:ODE_bad_set}
  If $u \in \Udeg$ then there exists a nonempty open subset $I_0 \subset I$ such that the restriction of $u$ to $I_0$ is continuously differentiable $r+1$ times, i.e., $\left. u \right\vert_{I_0} \in C^{r+1}(I_0)$. Furthermore, the restriction of $\mcal{L}(u)$ to $I_0$ is continuously differentiable.
\end{theorem}

\begin{proof}
  Let $ u \in \Udeg$, by definition there is a nonzero $\vect{a} \in \R^{N_1\times N_2}$ satisfying
  \begin{equation} \label{eq:ODE_degenerate}
    \sum_{i=1}^{N_1}\sum_{j=1}^{N_2} a_{i,j} \theta_{i}(t) \psi_{j}(\jet^r u(t)) = 0, \quad \text{for } t\in I.
  \end{equation}
  Gathering terms on each $\theta_{i}$ and denoting
  \[
    D(t, y_0, \ldots, y_r) = \sum_{i=1}^{N_1}\theta_{i}(t) \underbrace{\sum_{j=1}^{N_2} a_{i,j}  \psi_{j}(y_0, \ldots, y_r)}_{\phi_i(y_0, \ldots, y_r)},
  \]
  we have $
    D(t, \jet^r u(t)) = 0$ for all $t \in I$. We prove the theorem by induction on the maximum degree $d$ among the polynomials $\phi_{i} : \R^{r+1} \to \R$.
  We first observe that $d \geq 1$.
  If this were not the case, then all $\phi_{i}$ are constants equal to zero by linear-independence of the $\theta_{i}$.
  By linear independence of the $\psi_{j}$, this implies all coefficients $a_{i,j} = 0$, which is a contradiction.

  We begin with the $d=1$ case, meaning that every $\phi_{i}(y_0, \ldots, y_r) = b_{i} + \sum_{\ell=0}^r c_{i, \ell} y_\ell$ is an affine polynomial. That is,
  \[
    D(t,y) = \sum_{i=1}^{N_1} b_{i} \theta_{i}(t) + \sum_{\ell=0}^{r} \underbrace{\left( \sum_{i=1}^{N_1} c_{i, \ell} \theta_{i}(t) \right)}_{\partial_\ell D(t)} y_\ell,
    \quad \mbox{where} \quad \partial_\ell \coloneqq \frac{\partial}{\partial y_\ell}.
  \]
  If $\partial_\ell D = 0$ then $c_{i, \ell} = 0$ for all $1\leq i\leq N_1$ due to linear independence of the $\theta_{i}$. This implies that for each $1\leq i\leq N_1$,  $\phi_{i}$ does not depend on $y_\ell$.
  If $\partial_\ell D = 0$ for every $\ell$ then all $\phi_{i}$ have degree zero, which is a contradiction.
  Let $k$ be the largest integer satisfying $\partial_k D \neq 0$.
  By continuity of the functions $\theta_{i}$, there is a nonempty open subset $I_0 \subset I$ such that $\partial_k D(t) \neq 0$ for every $t \in I_0$.
  Moreover, using \cref{eq:ODE_degenerate} we have
  \[
    0 = D(t,\jet^r u(t)) = \sum_{i=1}^{N_1} b_{i} \theta_{i}(t) + \sum_{\ell=0}^{k} \partial_\ell D(t) \frac{\d^\ell}{\d t^\ell} u(t),
  \]
  since $\partial_\ell D=0$ for $\ell > k$.
  Therefore, for every $t\in I_0$ we have
  \[
    \frac{\d^k}{\d t^k} u(t) =
    -\frac{1}{\partial_k D(t)} \left[ \sum_{i=1}^{N_1} b_{i} \theta_{i}(t) + \sum_{\ell=0}^{k-1} \partial_\ell D(t) \frac{\d^\ell}{\d t^\ell} u(t) \right],
  \]
  where each $\partial_\ell D$ is smooth because every $\theta_{i}$ is smooth.
  Since $u \in C^r(I_0)$, the right-hand-side of the above expression defines a function in $C^{r + 1 - k}(I_0)$. Therefore, we have
  $\left. u \right\vert_{I_0} \in C^{r+1}(I_0)$,
  completing the proof of the $d=1$ case.

  We now assume that the theorem holds when the polynomials $\phi_{i}$ have degrees at most $d$.
  Suppose that the maximum degree among the polynomials $\phi_{i}$ is $d+1$, and let $k$ be the largest integer satisfying $
    \partial_k D = \sum_{i=1}^{N_1} \theta_{i} \partial_k \phi_{i}(y) \neq 0$.
  We observe that the polynomials $\partial_\ell \phi_{i}$ have degree at most $d$.
  If
  \[
    \partial_k D(t, \jet^r u(t))
    = \sum_{i=1}^{N_1}\sum_{j=1}^{N_2} a_{i,j} \theta_{i}(t) \partial_k \psi_{j} (\jet^r u(t))
    = 0, \quad \forall t \in I,
  \]
  then the conclusions of the theorem hold by the induction hypothesis.
  Therefore, without loss of generality we assume that there exists $t_0 \in I$ such that $\bar{y} = \jet^r u(t_0)$ satisfies $\partial_k D(t_0, \bar{y}) \neq 0$.
  For each $\ell > k$ we have $\partial_\ell D = 0$, meaning we can express
  $D(t, y) = \tilde{D}(t, y_{0:k})$, where $
    y_{0:k} \coloneqq (y_0, \ldots, y_k)$,
  and we have $\tilde{D}(t, \jet^k u(t)) = 0$ for all $t\in I$, and $\partial_k\tilde{D}(t_0, \bar{y}_{0:k}) \neq 0$,
  where $\bar{y}_{0:k} = \jet^k u(t_0)$. By the implicit function theorem there are open neighborhoods $U \ni (t_0, \bar{y}_{0:k-1})$ and $V \ni \bar{y}_k$, and a smooth function $h:U \to V$ so that the only solutions of $\tilde{D}(t,y_{0:k}) = 0$ with $(t,y_{0:k-1}, y_k) \in U \times V$ are given by $\tilde{D}(t, y_{0:k-1}, h(t,y_{0:k-1})) = 0$ for all $(t, y_{0:k-1}) \in U$. Thus, for all $t$ in the nonempty open subset
  $
    I_0 = \{ t \in I \mid (t, \jet^k u(t)) \in U \times V \}$
  we have
  \[
    \frac{\d^k}{\d t^k} u(t) = h(t, \jet^{k-1} u(t)).
  \]
  Since $h$ is smooth over $U$, the image of $\jet^{k-1} \left. u\right\vert_{I_0}$ lies in $U$, and $u \in C^r(I)$, it follows that the right-hand-side of the above expression defines a function in $C^{r+1-k}(I_0)$.
  Therefore, we have $\left. u \right\vert_{I_0} \in C^{r+1}(I_0)$, from which it is also evident that
  \[
    \left.\mcal{L}(u)\right\vert_{I_0} = \sum_{i=1}^{N_1}\sum_{j=1}^{N_2} c^*_{i,j} \left.\theta_{i}\right\vert_{I_0} \psi_{j}\circ \left.\jet^r u\right\vert_{I_0} \in C^{1}(I_0).
  \]
  The proof is complete by induction on $d$.
\end{proof}

The following corollary of \cref{thm:ODE_bad_set} provides deterministic ``rough'' inputs $f$ that do not satisfy the required regularity, and thus $f \notin \Fdeg(\vect{c}^*)$. Hence, a rough input $f \in C(I)$, in the sense we make precise below, ensures that the coefficients $\vect{c}^*$ can be identified using any corresponding solution $u \in C^r(I)$ with that input.

\begin{corollary}[Rough inputs are nondegenerate]\label{cor:rough_inputs}
  If $f \in C(I)$ is \emph{rough} in the sense that there is a dense set of points in $I$ where $f$ is non-differentiable or has a discontinuous derivative, then $f \notin \Fdeg(\vect{c}^*)$.
\end{corollary}

\begin{proof}
  Let $f\in C(I)$ be a rough function and suppose for contradiction that $f \in \Fdeg(\vect{c}^*)$. Then, there exists a solution $u \in C^r(I)$ of \cref{eqn:nonlinear_ODE} such that $L_u$ is not injective. By \cref{thm:ODE_bad_set}, there exists a nonempty open subset $I_0 \subset I$ such that $\left. u \right\vert_{I_0} \in C^{r+1}(I_0)$ and $\left. f \right\vert_{I_0} = \left. \mcal{L}(u) \right\vert_{I_0} \in C^1(I_0)$. Since $f$ is rough, there exists a point $x_0\in I_0$ such that $f$ is not differentiable or has a discontinuous derivative at $x_0$, hence $\left. f \right\vert_{I_0} \notin C^1(I_0)$, giving the required contradiction.
\end{proof}

For example, continuous nowhere differentiable functions such as the Weierstrass function do not lie in $\Fdeg(\vect{c}^*)$, and can therefore be used confidently as inputs to recover the coefficients in \eqref{eqn:nonlinear_ODE}.
Furthermore, it is shown in \cite{hunt1994prevalence} that continuous nowhere differentiable functions are \emph{prevalent} in the sense of \cite{hunt1992prevalence, christensen1972sets} among continuous functions on the unit interval $[0,1] \subset \R$, providing a prevalent set of recovering inputs when $I$ is bounded.

\subsection{Uniformly elliptic PDE with known higher-order term}

We now consider the problem of recovering the constant coefficients $\vect{b}^* \in \R^d$ and $c^*\in \R$ of a uniformly elliptic PDE in divergence form~\cite[Chapt.~6]{evans2022partial} on an open bounded domain $\Omega \subset \R^d$ with Lipschitz smooth boundary, where the higher-order variable coefficients matrix $\vect{A}\in L^\infty(\Omega)^{d\times d}$ is already known. That is, we consider the following partial differential equation:
\begin{equation}\label{eq:known_higher_order}
  \mathcal{L}u = -\nabla \cdot (\vect{A} \nabla u)+\vect{b}^* \cdot \nabla u + c^* u = f,\qquad u|_{\partial \Omega} = 0,
\end{equation}
where $f\in \mcal{F} = H^{-1}(\Omega)$ and $u\in \mcal{U} = H_0^1(\Omega)$. We assume that $\vect{b}^*$ and $c^*$ are chosen such that there is a unique weak solution $u\in H_0^1(\Omega)$ to~\cref{eq:known_higher_order} for each $f\in H^{-1}(\Omega)$ by the Lax--Milgram theorem~\cite[Sec.~6.2.1]{evans2022partial}, and the associated solution operator $\mcal{G}: \mcal{F} \to \mcal{U}$ is well-defined and continuous. Note that \cref{eq:known_higher_order} can be reformulated in the form of \cref{eqn:Lu_f_differential_ops} with $N=d+1$ as
\[\sum_{i=1}^d b_i^* \partial_i u + c^* u = \underbrace{f+\nabla \cdot (\vect{A} \nabla u)}_{\mathcal{D}_0(u,f)},\]
where $\mcal{D}_i(u,f) = \partial_i u$ and $\mcal{D}_{N}(u,f) = u$. The associated parameter-to-output recovery map $L_u: \R^d \times \R \to H^{-1}(\Omega)$ is defined as
\[
  L_u: (\vect{b}, c) \longmapsto \vect{b} \cdot \nabla u + c u,
\]
with the degenerate input set $\Fdeg(\vect{b}^*, c^*) = \{f \in H^{-1}(\Omega) \mid L_{\mcal{G}(f)} \text{ is not injective} \}$. According to \cref{cor:analytic_recovery}, the degenerate set is either compactly slice-null or equal to the entire space $\mathcal{F}$.  The next proposition constructs a solution $u_0 \in \mathcal{U}$ such that $L_{u_0}$ is injective, which implies that $\Fdeg(\vect{b}^*, c^*)$ is compactly slice-null.

\begin{proposition} [Existence of non-degenerate inputs]\label{prop:degenerate_set_empty}
  The degenerate input set $\Fdeg(\vect{b}^*, c^*)$ associated with \cref{eq:known_higher_order} is compactly slice-null and has measure zero with respect to any nondegenerate Gaussian measure on Borel sets of $\mcal{F}$.
\end{proposition}

\begin{proof}
  First note that since $H_0^1(\Omega)$ is a separable Hilbert space, its dual $H^{-1}$ is also a separable Hilbert space, hence admits a Schauder basis, and the assumptions of \cref{cor:analytic_recovery} are satisfied. To show that $\Fdeg(\vect{b}^*, c^*)\neq \mathcal{F}$, we construct a function $u_0\in H_0^1(\Omega)$ such that $L_{u_0}$ is injective, which implies that $f_0 = \mathcal{L}u_0 \notin \Fdeg(\vect{b}^*, c^*)$. Let $x_0\in\Omega$, since $\Omega$ is open, there exists a ball $B_\varepsilon(x_0)\subset\Omega$ of radius $\varepsilon>0$. Let $\phi\in C_c^\infty(\Omega)$ be a smooth bump function supported in $B_\varepsilon(x_0)$ such that $\phi=1$ on $B_{\varepsilon/2}(x_0)$, and consider the function $u_0:x\mapsto \phi(x)(x_1^2+\ldots+x_d^2)$ in $H_0^1(\Omega)$. Let $\vect{b}\in\R^d$, $c\in \R$ such that $L_{u_0}(\vect{b},c)=\vect{b}\cdot\nabla u_0 + c u_0 = 0$ in $\Omega$. After restricting to $B_{\varepsilon/2}(x_0)$, we have $L_{u_0}(\vect{b},c) = 2\vect{b}\cdot x + c \|x\|^2 = 0$ for all $x\in B_{\varepsilon/2}(x_0)$. This implies that $\vect{b}=0$ and $c=0$ since the function $2\vect{b}\cdot x + c \|x\|^2$ is a nontrivial polynomial in $x$ unless $\vect{b}=0$ and $c=0$.
\end{proof}

The next proposition shows that the degenerate input set $\Fdeg(\vect{b}^*, c^*)$ is actually trivial and that one can always recover the coefficients $\vect{b}^*$ and $c^*$ from any input-solution pair $(f,u) \in \mcal{F} \times \mcal{U}$ provided that $f \neq 0$ and $u \neq 0$.

\begin{proposition}[Trivial degenerate set]
  The degenerate input set associated with \cref{eq:known_higher_order} is trivial, i.e., $\Fdeg(\vect{b}^*, c^*) = \{0\}$.
\end{proposition}

\begin{proof}
  Let $f\neq 0\in H^{-1}(\Omega)$ and $u = \mathcal{G}(f) \neq 0 \in H_0^1(\Omega)$, we aim to show that the map $L_{\mathcal{G}(f)} = L_u$ is injective. Let $\vect{b}\in \R^d$ and $c \in \R$ such that $L_u(\vect{b},c) = \vect{b}\cdot\nabla u + c u = 0$ in $\Omega$. First note that if $\vect{b}=0$, then $cu = 0$, which implies that $c=0$ since $u\neq 0$, and the conclusion follows.

  Now suppose that $\vect{b}\neq 0$ by contradiction, since the equation $\vect{b}\cdot\nabla u + c u = 0$ is invariant under an orthogonal change of variables, we can assume without loss of generality that $\vect{b} = |\vect{b}|\vect{e}_1$. Hence, $|\vect{b}|\partial_{1} u + c u = 0$ in $\Omega$, or equivalently, $\partial_{1} u = -\tfrac{c}{|\vect{b}|} u$. Since $u \in H_0^1(\Omega)$, we extend $u$ by zero outside $\Omega$ to obtain $u \in H^1(\R^d) = W^{1,2}(\R^d)$, and note that $\vect{b}\cdot\nabla u + c u = 0$ holds in $L^2(\R^d)$. We now apply a Sobolev slicing theorem \cite[Thm.~4.21]{evans2015measure}, which states that for almost every $x' = (x_2,\ldots,x_d) \in \R^{d-1}$ with respect to the Lebesgue measure on $\R^{d-1}$, the slice $u_{x'}: t \mapsto u(t, x')$ is absolutely continuous on every compact subset of $\R$, and satisfies the ODE
  \[
    \frac{d}{dt} u_{x'}(t) = \partial_1 u(t, x') =-\frac{c}{|\vect{b}|} u_{x'}(t), \quad\text{for almost every }t \in \R.
  \]
  Note that the ODE holds trivially outside of the slice intersecting $\Omega$ since $u_{x'} \equiv 0$ there.
  By the fundamental theorem of calculus and since $u_{x'}$ is absolutely continuous, this integrates to $u_{x'}(t) = u_{x'}(0)\, e^{-\frac{c}{|\vect{b}|} t}$ for all $t \in \R$. Since $u \in L^2(\R^d)$, Fubini's theorem implies that $u_{x'} \in L^2(\R)$ for almost every $x'$. However, the function $t \mapsto u_{x'}(0)\,e^{-\frac{c}{|\vect{b}|}t}$ is in $L^2(\R)$ only if $u_{x'}(0) = 0$, which implies that $u_{x'} \equiv 0$ for almost every $x' \in \R^{d-1}$. Hence $u = 0$ in $\Omega$, which is a contradiction. Therefore $\vect{b} = 0$ and the conclusion follows.
\end{proof}

\subsection{Second-order elliptic PDEs with unknown coefficients} \label{sec:second_order_PDE}

In this section, we analyze the impact of unknown coefficients in a second-order elliptic equation on the structure of the degenerate input set. Let $a^*,b^*>0$ and consider the following equation
\begin{equation}\label{eq:2nd_order_elliptic}
  \mathcal L_{a^*, b^*} u =
  a^* u_{xx} + b^* u_{yy}=f, \quad x\in\Omega = (0,1)^2,
\end{equation}
with homogeneous boundary conditions, input function $f\in \mathcal{F} = L^2(\Omega)$, and solution $u\in \mathcal{U} =  H^2(\Omega) \cap H_0^1(\Omega)$. Here, the degenerate solution set $\Udeg$ is defined as the set of solutions $u\in  \mathcal{U}$ for which the parameter-to-observable map $L_u:(a,b)\mapsto a u_{xx} + b u_{yy}$, defined on $\R^2$, is not injective. The degenerate input set $\Fdeg(a^*, b^*)$ is defined as the set of inputs $f\in\mathcal{F}$ such that $L_{\mathcal{G}(f)}$ is not injective, where $\mathcal{G}:L^2(\Omega) \to H^2(\Omega) \cap H_0^1(\Omega)$ is the solution operator associated with \cref{eq:2nd_order_elliptic}. It follows from \cref{cor:analytic_diff_op_u_degenerate,cor:analytic_recovery_strong}  that both $\Udeg$ and  $\Fdeg(a^*, b^*) = \mathcal{L}_{a^*,b^*}(\Udeg)$ are compactly slice-null and have measure zero with respect to any nondegenerate Gaussian measure on Borel sets of $\mcal{U}$ and $\mcal{F}$, respectively. The next proposition explicitly characterizes the structure of the degenerate input and solution sets as  countable unions of closed subspaces.

\begin{proposition}[Structure of degenerate sets]\label{prop:diagonaL_2d}
  Consider the second-order elliptic operator $\mathcal{L}_{a^*,b^*}$ defined in \cref{eq:2nd_order_elliptic} for $a^\ast, b^\ast>0$. The degenerate solution and input sets $\Udeg$ and $\Fdeg(a^*, b^*)$ are given by
  \begin{align*}
    \Udeg           & = \!\!\!\!\! \bigcup_{\substack{m,n\geq 1\\ \gcd(m,n) = 1}}\!\!\!\!\overline{\vspan\{\psi_{km}\otimes \psi_{kn}\mid k\geq 1\}},\quad \psi_n:x\mapsto \sin(n\pi x), \\
    \Fdeg(a^*, b^*) & = \!\!\!\!\! \bigcup_{\substack{m,n\geq 1\\ \gcd(m,n) = 1}}\!\!\!\!\overline{\vspan\{\psi_{km}\otimes \psi_{kn}\mid k\geq 1\}}.
  \end{align*}
  Here, the closure is taken in $H^2(\Omega)$ for $\Udeg$ and in $L^2(\Omega)$ for $\Fdeg(a^*,b^*)$.
\end{proposition}

\begin{proof}
  We begin by applying the spectral theorem for compact self-adjoint operators to the solution operator associated with $u\mapsto -u_{xx}$ on $[0, 1]$ with homogeneous Dirichlet boundary conditions~\cite[Sec.~6.5]{evans2022partial}. The eigenfunctions of this operator are given by $\{\psi_i: x\mapsto\sin(i \pi x) \}_{i = 1}^\infty$, with corresponding eigenvalues $\{i^2 \pi^2\}_{i= 1}^\infty$, and form an orthogonal basis for $L^2([0, 1])$. Following~\cite[Sec.~II.4]{reed1980methods}, the set of pairwise products of the eigenfunctions $\{\psi_i \otimes \psi_j\}_{i,j=1}^\infty$ form an orthogonal basis for $L^2([0,1])\otimes L^2([0,1]) \cong L^2([0, 1]^2)$. Then, for all $i, j \geq 1$, the function $\psi_i \otimes \psi_j$ is an eigenfunction of $\mathcal L_{a, b}$ with eigenvalue $-(a i^2 \pi^2 + b j^2 \pi^2)$.

  Let $u\in \Udeg$. Then, the operator $L_u:(a,b)\mapsto a u_{xx}+b u_{yy}$ is not injective, hence there exists $(a,b)\neq (0,0)$ in the kernel of $L_u$, such that $a u_{xx}+b u_{yy} = 0$. After expanding $u$ in the basis $\{\psi_i\otimes \psi_j\}$ as $u = \sum_{i,j=1}^\infty u_{i,j}\psi_i\otimes \psi_j$, we obtain
  \[
    \mathcal{L}_{a,b}u = \sum_{i,j=1}^\infty -(a i^2 \pi^2 + b j^2 \pi^2) u_{i,j}\psi_i\otimes \psi_j= 0.
  \]
  Since the basis is orthogonal, we have $(a i^2 + b j^2)u_{i,j} = 0$ for all $i,j\geq 1$. Therefore, $u_{i,j} = 0$ if $a i^2 + b j^2 \neq 0$. This implies that $u\in \overline{\vspan\{\psi_{i}\otimes \psi_j\mid (i,j)\in\mathbb{Z}_{a,b}\}}$, where the closure is in $H^2(\Omega)$ and $\mathbb{Z}_{a,b} = \{(i,j)\in\mathbb{Z}_{\geq 1}^2\mid a i^2 + bj^2 = 0\}$. We note that each such closed subspace is contained in $H^2(\Omega)\cap H_0^1(\Omega)$ since $\psi_{i}\otimes \psi_j$ vanishes on $\partial [0,1]^2$ and $H^2(\Omega)\cap H_0^1(\Omega)$ is closed in $H^2(\Omega)$.

  Observe that  $\mathbb{Z}_{a,b}=\emptyset$ unless $(a,b)\neq (0,0)$ and $a/b\in \mathbb{Q}$ with $ab<0$. Then, $(i,j)\in \mathbb{Z}_{a,b}$ implies that $i^2/j^2 = -b/a$, hence $-b/a$ is the square of a rational number: $-b/a = (m/n)^2$, where $m, n
    \geq 1$ and $\gcd(m, n) = 1$. Therefore, $i/j = m/n$ and $(i,j) = (km,kn)$ for some $k\geq 1$ since $m,n$ are coprime, which implies that $\mathbb{Z}_{a,b} = \{(km,kn)\mid k\geq 1\}$. This shows that
  \[
    \Udeg\subset \!\!\!\!\!\bigcup_{\substack{m,n\geq 1\\ \gcd(m,n) = 1}}\!\!\!\!\overline{\vspan\{\psi_{km}\otimes \psi_{kn}\mid k\geq 1\}}.
  \]
  Conversely, suppose $u$ belongs to the right-hand side set of the above equation, there exist $m,n\geq 1$ such that $u\in \overline{\vspan\{\psi_{km}\otimes \psi_{kn}\mid k\geq 1\}}$. Let $(a,b) = (n^2,-m^2)$, then $a(km)^2+b(kn)^2 = 0$, hence
  \[L_u(a,b) = \sum_{k=1}^\infty -\pi^2(a(km)^2+b(kn)^2) u_{km,kn}\psi_{km}\otimes \psi_{kn}=0,\]
  which implies that $L_u$ is non-injective so $u\in \Udeg$.

  Let $f\in \Fdeg(a^*, b^*)$, then $L_{\mathcal{G}(f)}$ is non-injective, i.e., $u = \mathcal{G}(f)\in \Udeg$. Then, we have the following for some coprime $m,n\geq 1$,
  \[f = -\pi^2(a^*m^2+b^*n^2)\sum_{k=1}^\infty k^2 u_{km,kn}\psi_{km}\otimes \psi_{kn}\in \overline{\vspan\{\psi_{km}\otimes \psi_{kn}\mid k\geq 1\}},\]
  where $(k^2 u_{km,kn})_{k\geq 1}\in \ell^2(\mathbb{R})$ since $f\in L^2(\Omega)$. The reverse inclusion follows similarly.
\end{proof}

\cref{prop:diagonaL_2d} shows that the degenerate input set does not depend on the ground truth coefficients $a^*,b^*$. Hence, any input that does not belong to this set allows for the recovery of all operators of the form of \cref{eq:2nd_order_elliptic}. A similar version of \cref{prop:diagonaL_2d}  extends to higher dimensions, however we are only able to explicitly characterize $\Udeg$ in the $d=2$ case. In higher dimensions, the non-injectivity of the map $L_u$ leads to a system of Diophantine equations whose closed form solution is difficult to express analytically.

\section{Conclusions} \label{sec:conclusions}

We have shown that one-shot coefficient recovery for analytic linearly parameterized systems obeys a sharp zero-one law. Once a single input produces an injective recovery map, the set of inputs that fail to identify the coefficients is negligible in both a finite-dimensional slicing sense and under every nondegenerate Gaussian measure. This general principle applies to matrix recovery, analytic dynamical systems, and differential operators. Furthermore, in certain model ODE and PDE settings the degenerate sets can even be described more explicitly. Our results imply that finding non-degenerate inputs is easy when it is possible, so that attempting to recover coefficients using a random input either succeeds or provides an almost sure certificate that one-shot recovery is impossible. Several questions remain open, including stable recovery under noise, near-best approximation within the dictionary class, the design of finite test functionals that preserve infinite-dimensional identifiability in weak formulations, the value of adaptive versus not-adaptive input data, and one-shot recovery for nonlinearly-parameterized systems with partially-observed states.

%\section*{Acknowledgments}

\bibliographystyle{jabbrv_siam}
\bibliography{references}
\end{document}